\documentclass{amsart}

\usepackage{amssymb}
\usepackage{amsmath}
\usepackage{amsthm}
\usepackage{amscd}
\usepackage{changepage}
\usepackage{enumitem}
\usepackage{fancybox}
\usepackage{fancyhdr}
\usepackage{float}
\usepackage{hyperref}
\usepackage{marvosym}
\usepackage{mathrsfs}
\usepackage[pdftex]{graphicx}
\usepackage{setspace}
\usepackage{subfig}
\usepackage{theoremref}
\usepackage{verbatim}
\usepackage{xy}

\newcommand{\pdv}[2]{\dfrac{\partial #1}{\partial #2}}
\newcommand{\xx}{\times}

\newcommand{\R}{\mathbb{R}}
 
\newcommand{\br}[1]{\left( #1 \right)}
\newcommand{\brs}[1]{\left[ #1 \right]}
\newcommand{\norm}[1]{\left\Vert #1 \right\Vert}
\newcommand{\abs}[1]{\left\vert #1 \right\vert}
\newcommand{\lb}[0]{\left\lbrace}
\newcommand{\rb}[0]{\right\rbrace}
\newcommand\BoldSquare{%
  \setlength\fboxrule{1.1pt}\setlength\fboxsep{0pt}\fbox{\phantom{\rule{5pt}{5pt}}}}

\makeatletter
\def\blfootnote{\xdef\@thefnmark{}\@footnotetext}
\makeatother

\newtheorem{thm}{Theorem}[section]
\newtheorem{prop}{Proposition}[section]
\newtheorem{lem}{Lemma}[section]
\newtheorem{cor}{Corollary}[section]
\newtheorem{innercustomthm}{Theorem}
\newenvironment{customthm}[1]
  {\renewcommand\theinnercustomthm{#1}\innercustomthm}
  {\endinnercustomthm}
\newtheorem{innercustomprop}{Proposition}
\newenvironment{customprop}[1]
  {\renewcommand\theinnercustomprop{#1}\innercustomprop}
  {\endinnercustomprop}

\theoremstyle{definition}

\newtheorem{deff}{Definition}[section]

\title[Hardy-Littlewood Adapted to the Harmonic Oscillator]{A
  Hardy-Littlewood Maximal Operator Adapted to the Harmonic Oscillator}
\author{Julian Bailey}

\date{}

\begin{document}

\begin{abstract}
 This paper constructs a Hardy-Littlewood type maximal operator adapted to
the Schr\"{o}dinger operator $\mathcal{L} := -\Delta +
\abs{x}^{2}$ acting on $L^{2}(\R^{d})$. It achieves this through the use of the Gaussian grid $\Delta^{\gamma}_{0}$,
constructed in \cite{Maas2012} with the Ornstein-Uhlenbeck operator in
mind. At the scale of this grid, this maximal operator will resemble the classical
Hardy-Littlewood operator. At a larger scale, the cubes of the maximal
function are decomposed into cubes from $\Delta^{\gamma}_{0}$ and
weighted appropriately. Through this maximal function, a new class of
weights is defined, $A_{p}^{+}$, with the property that for any $w \in
A_{p}^{+}$ the heat
maximal operator associated with $\mathcal{L}$ is bounded from
$L^{p}(w)$ to itself. This class contains any other known class that possesses
this property. In particular, it is strictly larger than $A_{p}$.
\end{abstract}

  \maketitle

\section{Introduction and Preliminaries}
\label{sec:intro}

\blfootnote{\textit{Key words and phrases.}
  Hardy-Littlewood; weights; harmonic oscillator; heat maximal
  operator. \\
This research was partially supported by the Australian
Research Council through the Discovery Projects DP12010369 and
DP160100941.}

The Hardy-Littlewood operator is ubiquitous in classical harmonic
analysis. From the Lebesgue differentiation theorem to Calder\'{o}n-Zygmund theory, the importance of this averaging operator can hardly
be overstated. Classical harmonic analysis can be thought of as being
intricately linked to the Laplacian $\Delta$. Many of its fundamental
objects, including the Hardy-Littlewood operator, are closely related to the functional calculus of the
Laplacian. A current area of active research is the study of the harmonic
analysis associated with differential operators other than the
Laplacian. At present, there is no suitable candidate for the
Hardy-Littlewood operator in this setting. It is quite
possible that such an operator would play a fundamental role in
extending the theory even further. In this paper, our aim is the
construction of a Hardy-Littlewood
type maximal operator adapted to the Schr\"{o}dinger operator
$
\mathcal{L} := -\Delta + \abs{x}^{2}
$
on $L^{2}(\R^{d})$. In order to outline the details of this
construction, we must first present some motivating theory.

Note that throughout this paper, we will be working in the Euclidean
space $\R^{d}$ endowed with the Lebesgue measure $dx$. The dimension $d$ will
be considered to be fixed. Let $V : \R^{d} \rightarrow \R_{\geq 0}$ be a
potential that is non-identically zero and satisfies, for some $q > d
/ 2$ and $C > 0$, the reverse H\"{o}lder inequality,
$$
\br{\frac{1}{\abs{Q}} \int_{Q} V(y)^{q} dy}^{\frac{1}{q}} \leq
\frac{C}{\abs{Q}} \int_{Q} V(y) dy,
$$
for every cube $Q \subset \R^{d}$. Consider the Schr\"{o}dinger
operator $\mathcal{L}_{V} := - \Delta + V$ on $L^{2}(\R^{d})$.
An important step in the
comprehension of the harmonic analysis of such an operator was
made by Shen through  the introduction of the critical radius
function, see \cite{shen1995p}. This is defined by
$$
\rho_{V}(x) := \sup \lb r > 0 : \frac{1}{r^{d-2}} \int_{B(x,r)} V \leq 1 \rb
$$
for $x \in \R^{d}$; where $B(x,r)$ is the ball in $\R^{d}$, centered
at $x$ and of radius $r$. At a scale smaller than this critical radius, the operators
associated with $\mathcal{L}_{V}$ behave ``locally'' like their
classical counterparts for the Laplacian. This indicates that if we
are to construct a Hardy-Littlewood type maximal operator for
$\mathcal{L}$, then our construction should resemble the classical
Hardy-Littlewood operator at this local scale. What should it look like
at a larger scale? In order to answer this question, we must briefly
delve into some Gaussian harmonic analysis.

As is quite frequent in mathematics, when studying a particular
object, it can be fruitful to change perspective by
studying an isomorphic object in a different setting. Let $d \gamma
(x) := \pi^{-d/2} e^{-\abs{x}^{2}} dx$ denote the Gaussian measure on $\R^{d}$. Gaussian
harmonic analysis is
the study of the Ornstein-Uhlenbeck operator, $\mathcal{O} := -\Delta
+ 2 x \cdot \nabla$, on the space
$L^{2}(\gamma)$ and its associated harmonic analysis. Its relevance to
the study of $\mathcal{L}$ is that through the isometry $U :
L^{2}(dx) \rightarrow L^{2}(\gamma)$, defined by
$$
U f(x) := \pi^{-d /4} e^{-\frac{\abs{x}^{2}}{2}} f(x),
$$
for $f \in L^{2}(dx)$ and $x \in \R^{d}$, the operators $\mathcal{L}$ and $\mathcal{O}$ become, more-or-less,
similar. See \cite{abu2006hermite} for further details. This
similarity allows for the transfer of geometric ideas between the
Gaussian and the harmonic oscillator setting.

A measure $\mu$ on $\R^{d}$ is said to be doubling if there exists
some $C > 0$ such that
\begin{equation}
\label{eq:doubling}
\mu \br{B(x, 2 r)} \leq C \mu \br{B(x,r)},
\end{equation}
for all $x \in \R^{d}$ and $r > 0$. Many of the constructions from classical harmonic analysis directly
rely on the fact that the Lebesgue measure is doubling. A fundamental
obstruction in the development of Gaussian harmonic analysis is that,
due to the non-doubling nature of the Gaussian measure, many of these constructions do
not directly translate to the Gaussian setting.
In their seminal paper \cite{Mauceri2007}, Mauceri and Meda made a crucial step in
this development by transposing the critical
radius over to Gaussian harmonic analysis. They introduced their concept of
admissibility. 

Let us introduce $\rho$ as shorthand notation for the
critical radius function of $\mathcal{L}$, $\rho_{\abs{x}^{2}}$. It is not too difficult to
see that $\rho(x) = \min \lb 1 , 1 / \abs{x} \rb$. A ball $B(x,r)$
is then said to be admissible if $r \leq \rho(x)$. The collection of
all admissible balls in $\R^{d}$, $\mathcal{B}$, possesses the
desirable property that there exists some $C > 0$ such that the Gaussian measure satisfies the doubling
condition \eqref{eq:doubling} for all balls in $\mathcal{B}$. As such, by
restricting their attention to the collection $\mathcal{B}$, Mauceri
and Meda were
able to construct Gaussian analogues of the spaces $BMO$ and $H^{1}$.
A similar construction for the harmonic oscillator, also based on the
distinction between local and non-local scales, was developed by
Dziubanski and Zienkiewicz in \cite{dziubanski1999hardy} and
subsequent papers.

In \cite{Maas2012}, Maas, van Neerven and Portal extended the idea of admissibility by
constructing an admissible dyadic grid $\Delta^{\gamma}$. It is this grid that will form the
foundation for our construction. We recall some pertinent details.
For $m \in \mathbb{Z}$, let $\Delta_{m}$ denote the collection of
cubes
$$
\Delta_{m} := \lb 2^{-m} \br{x + [0,1)^{d}} : x \in \mathbb{Z}^{d} \rb.
$$
The standard dyadic grid is then the union $\Delta = \cup_{m \in
  \mathbb{Z}} \Delta_{m}$. Define the layers
$$
L_{0} := [-1,1)^{d}, \quad L_{l} := [-2^{l}, 2^{l})^{d} / [-2^{l-1}, 2^{l-1})^{d},
$$
for $l \geq 1$. Then define, for $k \in \mathbb{Z}$ and $l \geq 0$,
$$
\Delta^{\gamma}_{k,l} := \lb Q \in \Delta_{l+k} : Q \subseteq L_{l}
\rb, \quad \Delta^{\gamma}_{k} := \bigcup_{l \geq 0}
\Delta^{\gamma}_{k,l}, \quad \Delta^{\gamma} := \bigcup_{k \geq 0} \Delta^{\gamma}_{k}.
$$
The collection $\Delta^{\gamma}$ is called the Gaussian grid and will
be used extensively throughout this paper. Let's introduce some
notation that can be used
in conjunction with this grid. For any $x \in \R^{d}$, $R_{x}$ will be used
to denote the unique cube in
$\Delta^{\gamma}_{0}$ that contains the point $x$. For any $R \in
\Delta^{\gamma}_{0}$, $j(R)$ is defined to be the unique integer such
that $R \subset L_{j(R)}$. The more commonly used notation, $c_{Q}$ and
$l(Q)$, representing the center and side-length of a cube $Q$ respectively,
will also be used. Next we will define what will be considered to be
our local region in the Gaussian grid.

\begin{deff} 
 \label{def:NearRegion} 
 For a cube $R \in \Delta^{\gamma}_{0}$, fix a subcollection
 $\mathcal{N}(R) \subset \Delta^{\gamma}_{0}$ that satisfies the following
 two properties. 
\begin{enumerate}
\item[$\bullet$] $\mathcal{N}(R)$ contains all cubes
  $R' \in \Delta^{\gamma}_{0}$ satisfying
$$
d(R,R') < 2^{-j(R)},
$$
where
$
d(R,R') := \inf \lb \abs{x - y} : x \in R \ \mathrm{and} \ y \in R' \rb.
$
\item[$\bullet$] The region
$$
N(R) := \bigsqcup_{R' \in \mathcal{N}(R)} R',
$$
is a cube of sidelength $2^{2}l(R)$.
\end{enumerate}
The notation $\mathcal{F}(R) := \Delta^{\gamma}_{0} / \mathcal{N}(R)$
and $F(R) := \R^{d} / N(R)$ will also be employed.
 \end{deff}

It is obvious that such a subcollection must exist for each cube. There might even be more
than one such example. This, however,
is unimportant. What is important, is that we fix $\mathcal{N}(R)$ from
the outset. Examples of subcollections that satisfy these properties
are illustrated below.

\begin{figure}[H]
\centering
\includegraphics[width = 5.0in, height=1.5in]{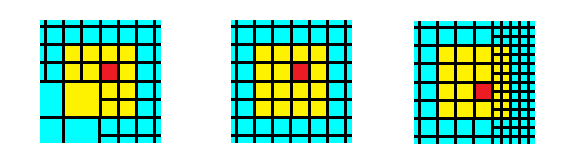}
\caption{Each of the above illustrations depicts a cube $R$, coloured
  in red, contained in the grid $\Delta^{\gamma}_{0}$ in dimension two. The near region, $N(R)$, consists of all
  cubes highlighted in yellow together with the cube $R$. The far region, $F(R)$, is coloured blue and extends out to infinity.}
\end{figure}

 As our operator is expected to behave differently at large
scales than at local scales, it is desirable to split it up into local
and non-local components. For any sub-linear operator
$B$, define
$$
B_{loc}f(x) := B \br{f \cdot \chi_{N(R_{x})}}(x) \quad \mathrm{and}
$$
$$
B_{far}f(x) :=  B \br{f \cdot \chi_{F(R_{x})}}(x), \quad \phantom{an}
$$
for $f \in L^{1}_{loc}(\R^{d})$ and $x \in \R^{d}$. Notice that
due to sub-linearity, for any weight $w$ on $\R^{d}$, to bound the quantity
$\norm{B f}_{L^{p}(w)}$, it is both sufficient and necessary to bound
$\norm{B_{loc}f}_{L^{p}(w)}$ and $\norm{B_{far}f}_{L^{p}(w)}$. 

Now that sufficient preliminaries have been discussed,
the details of our construction will be outlined. As noted previously,
$\Delta^{\gamma}_{0}$ acts as a mediator between the local and
non-local worlds. It is then appropriate to consider maximal functions
of the below general form as candidates for an adapted maximal
function for $\mathcal{L}$.

\begin{deff} 
 \label{def:GeneralMaximal}
For $Q \in \Delta$, let $\mathcal{G}(Q)$ be the collection of cubes
$$
\mathcal{G}(Q) := \lb \begin{array}{c c}
\lb Q \rb & \mathrm{if} \ Q \in \Delta^{\gamma}_{0} \\
\lb R' \in \Delta^{\gamma}_{0} : R' \subset Q \rb & \mathrm{otherwise}.
\end{array} \right.
$$
Then for $c : \Delta \xx \Delta^{\gamma}_{0} \xx \Delta^{\gamma}_{0}
\rightarrow \R_{\geq 0}$, $f \in L^{1}_{loc}(\R^{d})$ and $x \in R \in
\Delta^{\gamma}_{0}$, define the
operator $\mathcal{M}_{c}$ by
\begin{equation}
\label{eq:GeneralMaximal}
\mathcal{M}_{c}f(x) := \sup_{Q \in \Delta, \ Q \ni x} \frac{1}{\abs{Q}} \sum_{R'
  \in \mathcal{G}(Q)} c(Q,R,R') \int_{R'}
\abs{f(y)} dy.
\end{equation}
 \end{deff}

Notice that if $c(Q,R,R') = 1$ for all $R' \in \mathcal{G}(Q)$ and $Q
\in \Delta$, then the operator $\mathcal{M}_{c}$ is identical to the classical
dyadic Hardy-Littlewood operator. 

This looks promising but how do we determine what the right
$c$-coefficients are? Any candidate for an adapted Hardy-Littlewood
should share similar properties to the classic Hardy-Littlewood. We
will determine appropriate coefficients from one of these properties.
 Let $M$ and $T^{*}$ denote the classical
Hardy-Littlewood and heat maximal operator respectively. That is,
$$
M f(x) := \sup_{Q \ \mathrm{cube}, \ Q \ni x} \frac{1}{\abs{Q}} \int_{Q}
\abs{f(y)} dy \qquad \mathrm{and} \qquad T^{*}f(x) := \sup_{t > 0} e^{t \Delta} \abs{f}(x),
$$
for $f \in L^{1}_{loc}(\R^{d})$ and $x \in \R^{d}$. Also recall that
the $A_{p}$ class of weights is defined to be the collection of all
weights, $w$ on $\R^{d}$, for which there exists a constant $C > 0$ that
satisfies
$$
w(Q)^{\frac{1}{p}} \cdot w^{-\frac{1}{p-1}}(Q)^{\frac{p-1}{p}} \leq C \abs{Q}
$$
for all cubes $Q$ in $\R^{d}$.
The 
following theorem is a well-known result from weighted theory.
\begin{thm} 
  \label{thm:ClassicWeighted}
Let $w$ be a weight on $\R^{d}$ and $1 < p < \infty$. Then
$$
w \in A_{p} \quad \Leftrightarrow \quad \norm{M}_{L^{p}(w) \rightarrow
  L^{p}(w)} < \infty \quad \Leftrightarrow \quad
\norm{T^{*}}_{L^{p}(w) \rightarrow L^{p}(w)} < \infty.
$$
 \end{thm}

Refer to \cite{stein2016harmonic} sections V.4 and V.6 for proof. The above theorem indicates that if we are to construct a
Hardy-Littlewood type maximal operator for $\mathcal{L}$, then the
correct $c$-coefficients should satisfy the below equivalence for each
$1 < p < \infty$,
$$
\norm{\mathcal{T}^{*}}_{L^{p}(w) \rightarrow L^{p}(w)} < \infty \quad
\Leftrightarrow \quad \norm{\mathcal{M}_{c}}_{L^{p}(w) \rightarrow
  L^{p}(w)} < \infty,
$$
where $\mathcal{T}^{*}$ is the semigroup maximal operator associated to
$\mathcal{L}$,
$$
\mathcal{T}^{*}f(x) := \sup_{t > 0} e^{- t \mathcal{L}} \abs{f}(x).
$$
The coefficients for our generalised maximal function will be
optimised in an attempt to produce the above equivalence. 

A significant source of inspiration for this investigation stemmed
from \cite{bongioanni2011classes}. In this paper, Bongioanni, Harboure
and Salinas defined a new class of weights, $A_{p}^{\infty}$, for
which $\mathcal{T}^{*}$ was bounded on $L^{p}(w)$ for all weights $w
\in A_{p}^{\infty}$. What was interesting about this class was that
it was strictly larger than the classic Muckenhoupt class. It
seems that by including the potential $\abs{x}^{2}$, the weight class
$A_{p}$  effectively increases in size. It can be inferred from this,
that in order to produce a maximal function smaller than $M$ and
therefore a larger weight class, the
coefficients for our maximal function must be smaller than unity. The $A_{p}^{\infty}$ class is
defined to be $A_{p}^{\infty} := \cup_{\theta \geq 0} A_{p}^{\theta}$,
where $w \in A_{p}^{\theta}$ if
and only if there exists some constant $C > 0$ such that for all cubes
$Q \subset \R^{d}$,
$$
w(Q)^{1/p}w^{-\frac{1}{p-1}}(Q)^{\frac{p-1}{p}} \leq C \abs{Q} \br{1 +
\frac{l(Q)}{\rho(c_{Q})}}^{\theta}.
$$
In \cite{tang2015weighted}, the author developed a maximal function
$M^{\theta}$ adapted to the class $A_{p}^{\theta}$ in the sense that
$M^{\theta} : L^{p}(w) \rightarrow L^{p}(w)$ is bounded if and only if
$w \in A_{p}^{\theta}$. This operator is defined through
\begin{equation}
\label{eq:LRWMax}
M^{\theta}f(x) := \sup_{Q \ni x} \frac{1}{\psi_{\theta}(Q) \abs{Q}} \int_{Q} \abs{f(y)} dy,
\end{equation}
where
$$
\psi_{\theta}(Q) := \br{1 + \frac{l(Q)}{\rho(c_{Q})}}^{\theta}.
$$
Notice that this maximal function is also an example of the general
class from definition \ref{def:GeneralMaximal} with $c(Q,R,R') = \psi_{\theta}(Q)^{-1} < 1$. This allows for more weights in the class
$A_{p}^{\theta}$. However, it does not take into account the fact that if the cubes $R$ and $R'$ are
far apart, then the potential should have a larger effect and
therefore the coefficient $c(Q,R,R')$ should be smaller. The
coefficients that we define for our maximal function take this into account. 
 The main theorem of
this paper is stated below.

\begin{customthm}{A} 
 \label{thm:Main} 
 There exists maximal functions, $\mathcal{M}^{-}_{far}$ and
 $\mathcal{M}^{+}_{far}$, of similar form to definition
 \ref{def:GeneralMaximal} that satisfy the chain of implications
$$
\norm{\mathcal{M}^{+}_{far}}_{L^{p}(w)} < \infty
\quad \Rightarrow \quad \norm{\mathcal{T}^{*}_{far}}_{L^{p}(w)} < \infty
\quad \Rightarrow \quad \norm{\mathcal{M}^{-}_{far}}_{L^{p}(w)} < \infty,
$$
for any weight $w$ on $\R^{d}$ and $1 < p < \infty$.
 \end{customthm}

For a precise definition of the above maximal functions,
$\mathcal{M}^{-}_{far}$ and $\mathcal{M}^{+}_{far}$, and a proof of
this statement, refer to section \ref{sec:far}. A secondary result of
this paper that characterises the local behaviour of an adapted
maximal function is stated below.

\begin{customthm}{B} 
 \label{thm:Main2}
For any weight $w$ on $\R^{d}$ and $1 < p < \infty$,
$$
\norm{M_{loc}}_{L^{p}(w) \rightarrow L^{p}(w)} < \infty \quad
\Leftrightarrow \quad \norm{\mathcal{T}^{*}_{loc}}_{L^{p}(w) \rightarrow
  L^{p}(w)} < \infty.
$$ 
 \end{customthm}

This theorem will be proved in section \ref{sec:local}. Together, these
two statements demonstrate that for any weight in the class
$$
A_{p}^{+} := \lb w \ \mathrm{weight} \ \mathrm{on} \ \R^{d} : \norm{\mathcal{M}^{+}_{far}}_{L^{p}(w)
  \rightarrow L^{p}(w)} < \infty \ \mathrm{and} \
\norm{M_{loc}}_{L^{p}(w) \rightarrow L^{p}(w)} < \infty \rb,
$$
we have $\norm{\mathcal{T}^{*}}_{L^{p}(w) \rightarrow L^{p}(w)} <
\infty$. 

It is then natural to ask how our weight class compares with the class
$A_{p}^{\infty}$? Section \ref{sec:final} provides an answer to this
question in the form of the following proposition.

\begin{customprop}{C}
 \label{prop:Main3} 
 The following chain of strict inclusions holds for any $1 < p < \infty$,
$$
A_{p} \subsetneq A_{p}^{\infty} \subsetneq A_{p}^{+}.
$$
 \end{customprop}

The above inclusion indicates that our coefficients serve as an
improvement  upon the constant coefficients of \eqref{eq:LRWMax}.

Finally, in section \ref{sec:truncate}, the techniques developed throughout
this paper will be used to show that the heat maximal operator for
$\mathcal{L}$ can be safely truncated when considering weighted questions.

This paper is part of my PhD thesis, supervised by Pierre Portal at
the
Australian National University. It is inspired by discussions of my
supervisor with Paco Villarroya, aiming to understand better how to
adapt
harmonic analysis to the ¡°hidden geometry¡± of a differential
operator. In
most cases, this involves situations beyond the reach of
Calder\'{o}n-Zygmund
theory (see e.g.  \cite{Mauceri2007, hofmann2009hardy}). However, this can also be done within
Calder\'{o}n-Zygmund theory by proving stronger properties of smaller
classes of
singular integral operators than the Calder\'{o}n-Zygmund class. Paco
Villaroya
has particularly focused on describing compact (as opposed to merely
bounded)
singular integral operators (see e.g. \cite{villarroya2015characterization}). To do so, he had to
refine
classical dyadic approaches, in order to understand which cubes
particularly
affect compactness. I follow a similar path here, modifying standard
dyadic
arguments in a way that aims to reveal the ``hidden geometry'' of the
harmonic
oscillator. This is done by attempting to find the largest class of
weights
for which the corresponding heat maximal operator is bounded. Perhaps
surprisingly, such a question seems to be rarely formulated in the
context of
standard weighted Calder\'{o}n-Zygmund theory (generic questions involving
all
singular integrals and all related maximal functions are considered
instead),
but quite common in the context of two weights inequalities
(where studying just the Hilbert transform is hard enough, and
natural).

I would like to thank the anonymous referee of a previous
version of this paper for their suggestion that the proof of the
inclusion $A_{p}^{\infty} \subseteq A_{p}^{+}$ should be made strict.

\section{The Local Class}
\label{sec:local}

In this section, a local version of the $A_{p}$ class is introduced,
$A_{p}^{loc}$. This class is a dyadic variation of a similar class
introduced in \cite{bongioanni2011classes}. Through this class, and a
few preliminary lemmas,
Theorem \ref{thm:Main2} will be proved.

Consider a cube in $\R^{d}$, $Q_{0} := [a_{1}, a_{1} + l(Q_{0})) \xx
\cdots \xx [a_{d}, a_{d} + l(Q_{0}))$,
where $\lb a_{1}, \cdots, a_{d} \rb \subset \R$. In the usual manner, this cube can be divided into
$2^{d}$ congruent disjoint cubes with half the side-length of the original
cube. These cubes can themselves be divided into $2^{d}$ disjoint cubes
each and so on ad infinitum. If a cube $Q \subset \mathbb{R}^{d}$ can be obtained in
this manner from $Q_{0}$, then it is called a dyadic subcube of the cube
$Q_{0}$. Note that we did not require our initial cube $Q_{0}$ to be a member
of the standard dyadic grid and that $Q_{0}$ is a dyadic subcube of
itself.

\begin{deff} 
 \label{def:ApQ}
Fix a weight $w$ on $\mathbb{R}^{d}$ and $1 < p < \infty$. For a cube $Q_{0} \subset
\mathbb{R}^{d}$, the weight $w$ is said to belong to the class $A_{p}(Q_{0})$  if there exists a constant
$C > 0$ such that
\begin{equation}
\label{eq:near1}
w^{-\frac{1}{p-1}}(Q)^{\frac{p-1}{p}} w(Q)^{\frac{1}{p}} \leq C \abs{Q}
\end{equation}
for all dyadic subcubes $Q \subseteq Q_{0}$. The smallest such $C$ is
denoted $\brs{w}_{A_{p}(Q_{0})}$.
 \end{deff}

A variation of the next statement was originally proved in
\cite{harboure1984two}. It is an extension lemma for weights that
satisfy the $A_{p}$ property when restricted to a cube.

\begin{lem} 
 \label{lem:Extension}
Fix a cube $Q_{0} \subset \R^{d}$, $1 < p < \infty$ and a weight $w \in A_{p}(Q_{0})$. Then there exists a weight $w_{Q_{0}} \in
A_{p}(\R^{d})$ that coincides with $w$ on $Q_{0}$ such that
$\brs{w_{Q_{0}}}_{A_{p}} = \brs{w}_{A_{p}(Q_{0})}$.
\end{lem}

\begin{proof} 
Our proof proceeds by construction. Let $\mathcal{D}^{Q_{0}}$ denote a
dyadic system of cubes on $\R^{d}$ for which $Q_{0}$ is a
member. This can be explicitly constructed as follows. First, scale
the standard dyadic grid by a factor of $l(Q_{0})$ to form the
collection $l(Q_{0}) \cdot \Delta$ that consists of all cubes of the form
$$
 \left[ m_{1} 2^{k} l(Q_{0}), \br{m_{1} + 1} 2^{k} l(Q_{0})
  \right) \xx \cdots \xx \left[ m_{d} 2^{k} l(Q_{0}), \br{m_{d} + 1} 2^{k}
    l(Q_{0}) \right) 
  $$
  where $k, m_{1}, \cdots, m_{d} \in \mathbb{Z}$.
Then, if we let $b_{Q_{0}}$ denote the corner of the cube $Q_{0}$
closest to the origin, we can translate this scaled grid to $Q_{0}$,
$$
\mathcal{D}^{Q_{0}} := l(Q) \cdot \Delta + b_{Q_{0}} := \lb Q + b_{Q_{0}} :
Q \in l(Q) \cdot \Delta \rb.
$$
 Let $\mathcal{D}^{Q_{0}}_{0}$ denote the subcollection
that consists of all cubes in $\mathcal{D}^{Q_{0}}$ of the same size
as $Q_{0}$. A weight, $w_{Q_{0}}$ on $\R^{d}$, will be constructed for
which there exists $B > 0$ such that
\begin{equation}
\label{eq:Extension1}
w_{Q_{0}}^{-\frac{1}{p-1}}(Q)^{\frac{p-1}{p}}
w_{Q_{0}}(Q)^{\frac{1}{p}} \leq B \abs{Q}
\end{equation}
for all $Q \in \mathcal{D}^{Q_{0}}$. As the dyadic description of
$A_{p}(\R^{d})$ is scale and translation invariant, this criteria will
be sufficient to determine that $w_{Q_{0}} \in A_{p}(\R^{d})$. 

Fix $Q \in \mathcal{D}^{Q_{0}}_{0}$. Let $\varphi_{Q} : \R^{d}
\rightarrow \R^{d}$ denote the translation that takes the cube $Q$ to
the cube $Q_{0}$. Then, for $x \in Q$, define
$$
w_{Q_{0}}(x) := w \br{\varphi_{Q}(x)}.
$$
As the cubes in $\mathcal{D}^{Q_{0}}_{0}$ partition $\R^{d}$, this
description defines a unique function $w_{Q_{0}}$ on
$\R^{d}$. Moreover, it is clear that this function will be a weight
that coincides with $w$ on $Q_{0}$. 

By definition, as $w \in A_{p}(Q_{0})$, it follows that there must
exist a $C > 0$ such that \eqref{eq:Extension1} is satisfied for all
dyadic subcubes $Q \subset Q_{0}$. Fix a cube $Q \in
\mathcal{D}^{Q_{0}}$. Suppose that $Q$ is a dyadic subcube of a cube
from $\mathcal{D}^{Q_{0}}_{0}$. Then \eqref{eq:Extension1} must be
satisfied automatically with constant $C$. So suppose that $Q$ is not
a dyadic subcube of any cube in $\mathcal{D}^{Q_{0}}_{0}$. Then, since
a parent cube is always decomposable into its children, there must
exist finitely many cubes $\lb Q_{i} \rb_{i=1}^{N} \subset
\mathcal{D}^{Q_{0}}_{0}$ such that $Q = \sqcup_{i=1}^{N} Q_{i}$. We
then have
\begin{align*}\begin{split}  
 w_{Q_{0}}^{-\frac{1}{p-1}}(Q)^{\frac{p-1}{p}}
 w_{Q_{0}}(Q)^{\frac{1}{p}} &= \br{\int_{Q}
   w_{Q_{0}}(y)^{-\frac{1}{p-1}} dy}^{\frac{p-1}{p}} \br{\int_{Q}
   w_{Q_{0}}(y) dy}^{\frac{1}{p}} \\
&= \br{\sum_{i=1}^{N} \int_{Q_{i}} w_{Q_{0}}(y)^{-\frac{1}{p-1}}
  dy}^{\frac{p-1}{p}} \br{\sum_{i=1}^{N} \int_{Q_{i}} w_{Q_{0}}(y)
  dy}^{\frac{1}{p}} \\
&= \br{N \int_{Q_{0}} w^{-\frac{1}{p-1}}(y) dy}^{\frac{p-1}{p}} \br{N
  \int_{Q_{0}} w(y) dy}^{\frac{1}{p}} \\
&\leq C N \abs{Q_{0}} \\
&= C \abs{Q}.
 \end{split}\end{align*}
  \end{proof}

\begin{deff} 
 \label{def:ApLoc}
 Fix $1 < p < \infty$. A weight $w$ on $\R^{d}$ is said to be in the class $A_{p}^{loc}$ if
there exists a constant $C > 0$ such that
$$
\brs{w}_{A_{p}(N(R))} \leq C
$$
for all $R \in \Delta^{\gamma}_{0}$. The smallest such constant will
be denoted by $\brs{w}_{A_{p}^{loc}}$.
 \end{deff}

The subsequent lemma will be used numerous times throughout this
investigation. It states the exact form of the heat kernel
corresponding to $\mathcal{L}$. Its proof can be found in
\cite{simon1979functional} in dimension $1$. Higher dimensions follow
from this case by taking tensor products of Hermite functions.

\begin{lem} 
 \label{lem:Kernel} 
For $t > 0$, define the map $k_{t} : \R^{d} \xx \R^{d} \rightarrow \R$ through
\begin{equation}
\label{eq:kernel}
k_{t}(x,y) = h_{t}(x,y) \cdot \exp \br{- \alpha(t) \br{\abs{x}^{2} + \abs{y}^{2}}},
\end{equation}
where $h_{t}$ is the classic heat kernel
$$
h_{t}(x,y) := \frac{1}{\br{2 \pi t}^{d/2}} \exp \br{- \frac{\abs{x - y}^{2}}{2 t}}
$$
and $\alpha$ is defined by
$$
\alpha(t) := \frac{\sqrt{1 + t^{2}} - 1}{2 t}
$$
for all $x$ and $y$ in $\R^{d}$. The operator $\mathcal{T}^{*}$ is
then given by
$$
\mathcal{T}^{*}f(x) := \sup_{t > 0} \int_{\R^{d}} k_{t}(x,y)
\abs{f(y)} dy
$$
for any $f \in L^{1}_{loc}(\R^{d})$ and $x \in \R^{d}$. 
 \end{lem}

Note that the fundamental solution for $\mathcal{L}$ is actually
$k_{\sinh 2 t}$. We have chosen to rescale the kernel for simplicity. An expanded version of theorem \ref{thm:Main2} is presented and
proved below.

\begin{customthm}{B} 
 \label{thm:Local} 
 Let $T^{*}$ and $M$ denote the classic heat maximal operator and
 Hardy-Littlewood operator respectively. Let $w$ be a weight on $\R^{d}$. For any $1 < p < \infty$, the
 following statements are equivalent.
\begin{enumerate}
\item $\norm{M_{loc}}_{L^{p}(w) \rightarrow L^{p}(w)} < \infty$. \\
\item $w \in A_{p}^{loc}$. \\
\item $\norm{T^{*}_{loc}}_{L^{p}(w) \rightarrow L^{p}(w)} < \infty$. \\
\item $\norm{\mathcal{T}^{*}_{loc}}_{L^{p}(w) \rightarrow L^{p}(w)} < \infty$.
\end{enumerate}
\end{customthm}

\begin{proof}  
 We will prove the following chain of implications $(1) \Rightarrow
 (2) \Rightarrow (3) \Rightarrow (4) \Rightarrow (1)$. 

$(1) \Rightarrow (2)$.  Fix a cube $R \in \Delta^{\gamma}_{0}$, $Q$ a
dyadic subcube of $N(R)$ and $f \in L^{1}_{loc}(\R^{d})$. Define $C :=
\norm{M_{loc}}_{L^{p}(w) \rightarrow L^{p}(w)}$. Then,
using standard techniques from weighted theory,
\begin{align*}\begin{split}  
 \br{\int_{Q} w} \br{\frac{1}{\abs{Q}} \int_{Q} \abs{f}}^{p} &=
 \int_{Q} \br{\frac{1}{\abs{Q}} \int_{Q} \abs{f}}^{p} w(y) dy \\
&\leq \int_{Q} M_{loc}(f \cdot \chi_{Q})(y)^{p} w(y) dy \\
&\leq \norm{M_{loc}(f \cdot \chi_{Q})}^{p}_{L^{p}(\R^{n},w)} \\
&\leq C^{p} \norm{f \cdot \chi_{Q}}^{p}_{L^{p}(w)} \\
&= C^{p} \br{\int_{Q} \abs{f}^{p} w}.
 \end{split}\end{align*}
Take $f := \br{w + \varepsilon}^{-\frac{1}{p-1}}$ for some
$\varepsilon > 0$. Then
$$
w(Q) \br{\frac{1}{\abs{Q}} \int_{Q} \br{w(y) +
    \varepsilon}^{-\frac{1}{p-1}} dy}^{p} \leq C^{p} \int_{Q}
\frac{w(y)}{\br{w(y) + \varepsilon}^{\frac{p}{p-1}}} dy,
$$
which implies that
$$
w(Q) \br{\int_{Q} \br{w(y) + \varepsilon}^{-\frac{1}{p-1}} dy}^{p}
\leq C^{p} \abs{Q}^{p} \int_{Q} \frac{\br{w(y) + \varepsilon}}{\br{w(y) +
    \varepsilon}^{\frac{p}{p-1}} } dy,
$$
$$
\Rightarrow w(Q) \br{\int_{Q} \br{w(y) + \varepsilon}^{-\frac{1}{p-1}}
dy}^{p-1} \leq C^{p} \abs{Q}^{p}
$$
for each $\varepsilon > 0$. An application of the Lebesgue monotone
convergence theorem then produces the desired result. 

$(2) \Rightarrow (3)$. Lemma \ref{lem:Extension} states that for any
cube $R \in \Delta^{\gamma}_{0}$ the restriction $w \vert_{N(R)}$ can be
extended to an $A_{p}$ weight $w_{N(R)}$. As $w_{N(R)} \in A_{p}$, we know
from classical theory that $\norm{T^{*}}_{L^{p}(w_{N(R)}) \rightarrow
  L^{p}(w_{N(R)})} \lesssim \brs{w_{N(R)}}_{A_{p}} < \infty$. Then,
for $f \in L^{p}(w)$,
\begin{align*}\begin{split}  
 \norm{T_{loc}^{*}f}^{p}_{L^{p}(w)} &= \int_{\R^{d}}
 T_{loc}^{*}f(x)^{p} w(x) dx \\
&= \sum_{R \in \Delta^{\gamma}_{0}} \int_{R} T^{*}(f \cdot
\chi_{N(R)})(x)^{p} w(x) dx \\
&\leq \sum_{R \in \Delta^{\gamma}_{0}} \int_{\R^{d}} T^{*}(f \cdot
\chi_{N(R)})(x)^{p} w_{N(R)}(x) dx \\
&\lesssim \sum_{R \in \Delta^{\gamma}_{0}} \brs{w_{N(R)}}_{A_{p}}^{p} \int_{N(R)} \abs{f(x)}^{p}
w_{N(R)}(x) dx \\
&\leq \brs{w}_{A_{p}^{loc}}^{p} \sum_{R \in \Delta^{\gamma}_{0}} \int_{N(R)} \abs{f(x)}^{p} w(x) dx \\
&\lesssim \brs{w}_{A_{p}^{loc}}^{p} \int_{\R^{d}} \abs{f(x)}^{p} w(x) dx,
 \end{split}\end{align*}
where the final inequality was obtained from the bounded overlap
property of the cubes $\lb N(R) \rb_{R \in \Delta^{\gamma}_{0}}$.

$(3) \Rightarrow (4)$. This follows trivially from the inequality
$k_{t}(x,y) \leq h_{t}(x,y)$ for all $x , y \in \R^{d}$ and $t > 0$.

$(4) \Rightarrow (1)$. Fix $f \in L^{1}_{loc}(\R^{d})$ and $x \in
R \in \Delta^{\gamma}_{0}$. Let $Q$ be any cube containing $x$ that satisfies $Q \subseteq
N(R)$. We first observe that for any $y
 \in Q$,
$$
\exp(-\frac{\abs{x-y}^{2}}{2 l(Q)^{2}}) \approx 1.
$$
To see this, note that 
$$
\abs{x - y} \leq \sqrt{d} l(Q).
$$
This implies that
$$
- \frac{\abs{x-y}^{2}}{2 l(Q)^{2}} \geq -\frac{d}{2},
$$
and therefore
$$
\exp(- \frac{\abs{x- y}^{2}}{2 l(Q)^{2}}) \gtrsim 1.
$$
Moreover, we trivially have
$$
\exp(- \frac{\abs{x - y}^{2}}{2 l(Q)^{2}}) \leq 1.
$$
Note that for any $x , y \in Q$, since $l(Q) \leq 4 l(R)$, we have the
bound
\begin{align*}\begin{split}  
 \abs{x}, \abs{y} &\leq 2 \sqrt{d} 2^{j(R)} \\
&= \frac{2 \sqrt{d}}{l(R)} \\
&\leq \frac{8 \sqrt{d}}{l(Q)}.
 \end{split}\end{align*}
 This then implies that
\begin{align*}\begin{split}  
 \exp \br{-\frac{\br{\sqrt{1 + l(Q)^{4}}-1}}{2 l(Q)^{2}} \br{\abs{x}^{2} +
     \abs{y}^{2}}} \geq \exp \br{- \frac{8^{2} d \br{\sqrt{1 + l(Q)^{4}} - 1}}{
     l(Q)^{4}} }.
 \end{split}\end{align*}
It is easy to show that the bound
$$
\frac{\sqrt{1 + t^{4}} - 1}{t^{4}} \leq \frac{1}{2}
$$
is satisfied for all $t > 0$. This then gives us
$$
\exp \br{- \frac{\br{\sqrt{1 + l(Q)^{4}}-1}}{2 l(Q)^{2}} \br{\abs{x}^{2} +
    \abs{y}^{2}}} \geq e^{-\frac{8^{2} d}{2}}.
$$
For $t := l(Q)^{2}$, we then have
\begin{align*}\begin{split}  
 \frac{1}{\abs{Q}} \int_{Q}& \abs{f(y)} dy \\ &\lesssim \frac{1}{l(Q)^{d}}
 \int_{Q} \exp \br{- \frac{\br{\sqrt{1 + l(Q)^{4}} - 1}}{2 l(Q)^{2}}
   \br{\abs{x}^{2} + \abs{y}^{2}}}
 \exp \br{-\frac{\abs{x-y}^{2}}{2 l(Q)^{2}}} \abs{f(y)} dy \\
&= \int_{Q} \frac{1}{t^{d/2}} \exp \br{- \frac{\br{\sqrt{1 + t^{2}} - 1}}{2 t}
   \br{\abs{x}^{2} + \abs{y}^{2}}} \exp \br{-\frac{\abs{x - y}^{2}}{2 t}}
\abs{f(y)} dy \\
&\lesssim \int_{Q} k_{t}(x,y) \abs{f(y)} dy \\
&\lesssim \mathcal{T}^{*}_{loc}f(x).
 \end{split}\end{align*}
On taking the supremum over all such $Q$, we obtain $M_{loc}f(x)
\lesssim \mathcal{T}^{*}_{loc}f(x)$.
 \end{proof}

\section{The Far Class}
\label{sec:far}

In this section, the adapted operators $\mathcal{M}^{-}_{far}$ and
$\mathcal{M}^{+}_{far}$ are defined and theorem \ref{thm:Main} is proved. With
this, a sufficient condition  for the boundedness of
$\norm{\mathcal{T}^{*}}_{L^{p}(w) \rightarrow L^{p}(w)}$ is
obtained. Prior to presenting these definitions, it is necessary to
introduce a collection of cubes that represent the regions over which
our averaging operators will act.

\begin{deff} 
 \label{def:QtR}
For each $R \in \Delta^{\gamma}_{0}$, define the following subsets of
$\R^{d}$.
\begin{enumerate}
\item[$\bullet$] $Q_{0}(R)$ is the smallest cube containing the
  region
$$
\lb y \in \R^{d} : \abs{y} \leq 2^{16} d^{4} 2^{j(R)} \rb,
$$
that can be decomposed into cubes from the grid $\Delta^{\gamma}_{0}$.
\item[$\bullet$] For $t \leq 2^{4} d^{2}$, $Q_{t}(R) := Q_{0}(R)$.
\item[$\bullet$] For $t \geq 2^{4} d^{2}$, $Q_{t}(R)$ is the smallest
  cube containing the region
$$
\lb y \in \R^{d} : \abs{y} \leq 2^{8} t^{2} 2^{j(R)} \rb,
$$
that can be decomposed into cubes from the grid $\Delta^{\gamma}_{0}$.
\end{enumerate}
 \end{deff}

For sets $A$ and $B$ contained in $\R^{d}$, introduce the notation
$k_{t}^{+}(A,B)$ and $k_{t}^{-}(A,B)$ to denote respectively the
supremum and infimum of $k_{t}(x,y)$ over all $x \in A$ and $y \in B$.

\begin{deff} 
 \label{def:AdaptedMaximalFar} 
 For $f \in L^{1}_{loc}(\R^{d})$ and $x \in R \in
 \Delta^{\gamma}_{0}$, define the operators $\mathcal{M}^{+}_{far}$
 and $\mathcal{M}^{-}_{far}$ through

\vspace*{-0.15in}

\begin{align}\begin{split}
\label{eq:AdaptedMaximalFar}  
 &\mathcal{M}^{+}_{far}f(x) := \sup_{t > 0} \sum_{R' \in
  \mathcal{F}(R), \ R' \subset Q_{t}(R)} k_{t}^{+}(R,R')
\int_{R'} \abs{f(y)} dy, \qquad \mathrm{and} \\
\phantom{c} \\
&\mathcal{M}^{-}_{far}f(x) := \sup_{t > 0} \sum_{R' \in
  \mathcal{F}(R), \ R' \subset Q_{t}(R)} k_{t}^{-}(R,R')
\int_{R'} \abs{f(y)} dy.
 \end{split}\end{align}
 \end{deff}

With the introduction of our maximal functions, it is a
straightforward matter to define their corresponding weight classes.

\begin{deff} 
 \label{def:ApFar} 
 For $1 < p < \infty$, the classes of weights on $\R^{d}$, $A_{p}^{far+}$ and $A_{p}^{far-}$, are defined
 through

\vspace*{-0.15in}

\begin{align}\begin{split} 
\label{eq:ApPlus} 
 &A_{p}^{far+} := \lb w \ \mathrm{weight} \ \mathrm{on} \ \R^{d} : \norm{\mathcal{M}^{+}_{far}}_{L^{p}(w)
  \rightarrow L^{p}(w)} < \infty \rb  \quad \mathrm{and} \\ 
\phantom{c} \\
& A_{p}^{far-} := \lb w \ \mathrm{weight} \ \mathrm{on} \ \R^{d} : \norm{\mathcal{M}^{-}_{far}}_{L^{p}(w)
  \rightarrow L^{p}(w)} < \infty \rb. 
 \end{split}\end{align}

We then define $A_{p}^{+} := A_{p}^{far+} \cap A_{p}^{loc}$ and
$A_{p}^{-} := A_{p}^{far-} \cap A_{p}^{loc}$.
 \end{deff}

In order to verify our main result, a string of technical lemmas
must first be proved. The first two of these provide some valuable estimates
concerning the maximum of the function $t \mapsto k_{t}(x,y)$ for
fixed $x$ and $y$ in $\R^{d}$.

\begin{lem} 
 \label{lem:MaxKernel} 
 Fix points $ x \in R \in \Delta^{\gamma}_{0}$ and $y
 \notin Q_{0}(R)$. There is precisely one maximum for the function $t
 \mapsto k_{t}(x,y)$. Denote this point by $t_{m}(x,y)$. Then for $R$
 not contained in the first layer, $t_{m}(x,y)$ must satisfy
\begin{equation}
\frac{\abs{y}}{9 \cdot d \abs{x}} \leq t_{m}(x,y) \leq \frac{\abs{x -
    y}^{2}}{d}. 
\end{equation}
For $R$ contained in the first layer, $t_{m}(x,y)$ will satisfy
\begin{equation}
\frac{\abs{y}}{9 \cdot d} \leq t_{m}(x,y) \leq \frac{\abs{x - y}^{2}}{d}.
\end{equation}
\end{lem}

\begin{proof}  
On differentiating expression \eqref{eq:kernel} with respect to $t$ we obtain
$$
\pdv{}{t}k_{t}(x,y) = \frac{1}{2 t^{2}} g(t) k_{t}(x,y),
$$
where the function $g$ is defined to be
$$
g(t) := - d \cdot t + \frac{\br{\abs{x}^{2} + \abs{y}^{2}}}{\sqrt{1 +
    t^{2}}} - 2 \langle x , y \rangle.
$$
As the kernel $k_{t}(x,y)$ is always positive, it follows that the
sign of the derivative will be identical to the sign of the function
$g(t)$. Suppose that $g$ is negative. Then we must have
$$
 \br{d \cdot t +
2 \langle x , y \rangle} \sqrt{1 + t^{2}} > \br{\abs{x}^{2} + \abs{y}^{2}}
$$
That is, the derivative of the kernel will be negative if and only if
the above inequality holds. Likewise, the derivative of the kernel will
be positive if and only if 
\begin{equation}
\label{eq:Less}
 \br{d \cdot t + 2 \langle x , y \rangle} \sqrt{1 + t^{2}} <
\br{\abs{x}^{2} + \abs{y}^{2}}
\end{equation}
and the derivative will vanish if and only if equality holds.

It is simple to show that $\abs{x - y}^{2} / d$ serves as the only maximum of
the function $t \mapsto h_{t}(x,y)$. This implies that $h_{t}(x,y)$ is decreasing for $t > \abs{x - y}^{2} / d$. As the function
$\alpha(t)$ is strictly increasing, we have that 
$$
\exp \br{- \alpha(t)\br{\abs{x}^{2} + \abs{y}^{2}}}
$$
is strictly decreasing for all $t$. This shows that $k_{t}(x,y)$
is strictly decreasing for $t > \abs{x - y}^{2} / d$. It then
follows that any maximum for $t \mapsto k_{t}(x,y)$ must be less than
$\abs{x - y}^{2} / d$. As this function must approach
$0$ as $t$ approaches $0$, continuity of the derivative then
implies that there must exist at least one maximum in the interval
$\brs{0, \abs{x - y}^{2} / d}$.

 Let $t_{m}(x,y)$ denote the largest maximum in the above
interval. It will be shown that $t_{m}(x,y)$ is the only maximum. From our previous argument, equality will hold in  \eqref{eq:Less} for the value $t_{m}(x,y)$. Suppose that $t_{0} <
t_{m}(x,y)$. Then $t_{0} = t_{m}(x,y) - a$ for some $a > 0$. We then have
\begin{align*}\begin{split}  
 \br{d \cdot t_{0} + 2 \langle x , y \rangle} \sqrt{1 + t_{0}^{2}} &= \br{d
   \cdot t_{m}(x,y) - d \cdot a + 2 \langle x , y \rangle } \sqrt{1 +
   t_{0}^{2}} \\
&= \br{d \cdot t_{m}(x,y) + 2 \langle x , y \rangle} \sqrt{1 + t_{0}^{2}}
- d \cdot a \sqrt{1 + t_{0}^{2}}.
 \end{split}\end{align*}
As equality holds in expression \eqref{eq:Less} for $t_{m}(x,y)$, it
follows that the factor  $\br{d
  \cdot t_{m}(x,y)} + 2 \langle x , y \rangle$ must be
positive. Therefore
\begin{align*}\begin{split}  
 \br{d \cdot t_{0} + 2 \langle x , y \rangle} \sqrt{1 + t_{0}^{2}}
 &\leq \br{d \cdot t_{m}(x,y) + 2 \langle x , y \rangle} \sqrt{1 +
   t_{m}(x,y)^{2}} - d \cdot a \sqrt{1 + t_{0}^{2}} \\
&= \br{\abs{x}^{2} + \abs{y}^{2}} - d \cdot a \sqrt{1 + t_{0}^{2}} \\
&< \br{\abs{x}^{2} + \abs{y}^{2}}.
 \end{split}\end{align*}
This demonstrates that the derivative must be positive for any $t_{0} <
t_{m}(x, y)$.

Let's now show the lower bound for $t_{m}(x,y)$.
First suppose that $R$ is not contained in the first layer. It will be
shown that for any $t_{1} < \abs{y} / \br{ 9 \cdot d \abs{x}}$, inequality \eqref{eq:Less} holds. From our previous
argument, this will then imply that the function is increasing on the
interval $\left[0, \abs{y} / \br{ 9 \cdot d \abs{x}}\right)$.
 As $y \notin Q_{0}(R)$, it follows that $y$
satisfies the bound $\abs{y} > 3 \abs{x}$.  We know that
\begin{align*}\begin{split}  
 1 + t_{1}^{2} &<  1 + \frac{1}{9} \br{\frac{\abs{y}}{3 \abs{x}}}^{2} \\
 &=  1 +  \br{\frac{\abs{y}}{3   \abs{x}}}^{2} - \frac{8}{9} \br{\frac{\abs{y}}{3 \abs{x}}}^{2} \\
&\leq  1 +  \br{\frac{\abs{y}}{3 \abs{x}}}^{2} - \frac{8}{9} \\
&= \frac{1}{9} \br{1 + \frac{ \abs{y}^{2}}{\abs{x}^{2}} }.
 \end{split}\end{align*}
We also have
\begin{align*}\begin{split}  
 \br{ d \cdot t_{1} + 2 \langle x , y \rangle} &\leq \br{d \cdot t_{1}
 + 2 \abs{\langle x , y \rangle}} \\
&\leq \br{\frac{\abs{y}}{9 \abs{x}} +
   2 \abs{x} \abs{y}} \\
&\leq \br{\frac{\abs{y}}{\abs{x}}} \br{\frac{1}{9} + 2
  \abs{x}^{2}} \\
&\leq \br{\frac{\abs{y}}{\abs{x}}} 3 \abs{x}^{2}.
 \end{split}\end{align*}
This demonstrates that
\begin{align*}\begin{split}  
  \br{d \cdot t_{1} + 2 \langle x , y \rangle} \sqrt{1 + t_{1}^{2}} &<
\br{3 \abs{x} \abs{y}} \cdot \frac{1}{3} \sqrt{1 +
  \frac{\abs{y}^{2}}{\abs{x}^{2}}} \\
&= \abs{y} \sqrt{\abs{x}^{2} + \abs{y}^{2}} \\
&\leq \br{\abs{x}^{2} + \abs{y}^{2}}.
 \end{split}\end{align*}

Now suppose that $R$ is in the first layer and $y \notin
Q_{0}(R)$. Then $\abs{y} \geq 2^{16} d^{4}$. Let $t_{2} < \abs{y} / \br{
9 d}$. Then
\begin{align*}\begin{split}  
 \br{1 + t_{2}^{2}} &< \br{1 + \br{\frac{\abs{y}}{9 d}}^{2}} \\
&\leq \br{\frac{\abs{y}^{2}}{2^{32} d^{8}} + \frac{\abs{y}^{2}}{9^{2}
    d^{2}}} \\
&\leq \frac{2 \abs{y}^{2}}{9^{2} d^{2}}.
 \end{split}\end{align*}
On noting that
$\abs{x} \leq \sqrt{d}$,
\begin{align*}\begin{split}  
 \br{d \cdot t_{2} + 2 \langle x , y \rangle} &\leq \br{d \cdot t_{2}
   + 2 \abs{\langle x , y \rangle}} \\
&\leq
 \br{\frac{\abs{y}}{9} + 2 \abs{x} \abs{y}} \\
&\leq \br{\frac{1}{9} + 2 \abs{x}} \abs{y} \\
&\leq \br{\frac{1}{9} + 2 \sqrt{d}} \abs{y} \\
&\leq 3 \sqrt{d} \abs{y}
 \end{split}\end{align*}
This finally leads to
\begin{align*}\begin{split}  
  \br{d \cdot t_{2} + 2 \langle x , y \rangle} \sqrt{1 + t_{2}^{2}}
 &< \br{3 \sqrt{d} \abs{y}} \br{\frac{\sqrt{2} \abs{y}}{9 d}} \\
&\leq \abs{y}^{2} \\
&\leq \br{\abs{x}^{2} + \abs{y}^{2}},
 \end{split}\end{align*}
which validates our lower bound.
\end{proof}

\begin{lem} 
 \label{lem:Taylor} 
 Fix cubes $R$ and $R'$ in $\Delta^{\gamma}_{0}$ with $R' \subset
 Q_{0}(R)^{c}$. Fix points $x \in R$ and $y \in R'$. The maximum $t_{m}(x,y)$ satisfies the inequality,
\begin{equation}
2 \leq 8 \cdot t_{m}(x,y) \sqrt{\frac{2^{j(R) + j(R')}}{\abs{x}^{2} +
    \abs{y}^{2}}} \leq \frac{t_{m}(x,y)}{2^{4} d^{2}}.
\end{equation}
\end{lem}

\begin{proof}  
As $y \notin Q_{0}(R)$, we have $\abs{y} \geq 2^{16} d^{4}
2^{j(R)}$ and also $\abs{y} \geq 2^{j(R') - 1}$. The upper inequality
then follows from
\begin{align*}\begin{split}  
\abs{x}^{2} + \abs{y}^{2} &\geq \abs{y}^{2} \\
&\geq 2^{j(R')-1} 2^{16} d^{4} 2^{j(R)} \\
&= 2^{15} d^{4} 2^{j(R) + j(R')}.
 \end{split}\end{align*}
 As for the lower bound, first consider when $R$ is not in the first
layer. On applying Lemma \ref{lem:MaxKernel} and recalling
that $\abs{y} \geq \abs{x}$,
\begin{align*}\begin{split}  
 t_{m}(x,y) \sqrt{\frac{2^{j(R) + j(R')}}{\abs{x}^{2} + \abs{y}^{2}}}
 &\geq \frac{\abs{y}}{9 d \abs{x}} \sqrt{\frac{2^{j(R) +
       j(R')}}{\abs{x}^{2} + \abs{y}^{2}}} \\
&\geq \frac{1}{9 d} \sqrt{\frac{\abs{y}^{2} 2^{j(R)+j(R')}}{2
    \abs{x}^{2} \abs{y}^{2}}}.
 \end{split}\end{align*}
Then on applying the bounds $\abs{x} \leq \sqrt{d} 2^{j(R)}$, $\abs{y} \leq \sqrt{d}
2^{j(R')}$ and $\abs{y} \geq 2^{16} d^{4} 2^{j(R)}$ successively we
obtain
\begin{align*}\begin{split}  
 t_{m}(x,y) \sqrt{\frac{2^{j(R) + j(R')}}{\abs{x}^{2} + \abs{y}^{2}}} &\geq \frac{1}{9 d} \sqrt{\frac{2^{j(R) + j(R')}}{2 d 2^{2 j(R)}}}  \\
&\geq \frac{1}{9 d} \sqrt{\frac{\abs{y}}{2 d^{3/2} 2^{j(R)}}} \\
&\geq \frac{1}{9 d} \sqrt{\frac{2^{16} d^{4} 2^{j(R)}}{2 d^{3/2} 2^{j(R)}}} \\
&\geq 2.
 \end{split}\end{align*}
Next, consider when $R$ is in the first layer. Once again apply Lemma
\ref{lem:MaxKernel} and $\abs{y} \geq \abs{x}$ to obtain
\begin{align*}\begin{split}  
    t_{m}(x,y) \sqrt{\frac{2^{j(R) + j(R')}}{\abs{x}^{2} +
        \abs{y}^{2}}} &\geq \frac{\abs{y}}{9 d} \sqrt{\frac{2^{j(R')}}{2
      \abs{y}^{2}}} \\
&= \frac{1}{9 d} \sqrt{\frac{2^{j(R')}}{2}}.
 \end{split}\end{align*}
Then, on successively applying the bounds $\abs{y} \leq \sqrt{d}
2^{j(R')}$ and $\abs{y} \geq 2^{16} d^{4}$,
\begin{align*}\begin{split}  
 t_{m}(x,y) \sqrt{\frac{2^{j(R)} + 2^{j(R')}}{\abs{x}^{2} +
     \abs{y}^{2}}}
&\geq \frac{1}{9 d} \sqrt{\frac{\abs{y}}{2 \sqrt{d}}} \\
&\geq \frac{1}{9 d} \sqrt{\frac{2^{16} d^{4}}{2 \sqrt{d}}} \\
&\geq 2.
 \end{split}\end{align*}
This concludes the proof.
 \end{proof}

The next lemma obtains an estimate on ratios of the form $k_{t}(x,y) \cdot
k_{t_{m}(x,y)}(x,y)^{-1}$ for fixed $x$ and $y$. It will play a key role in
the proof of theorem \ref{thm:Main}.

\begin{lem} 
  \label{lem:FarKernel}
Fix cubes $R$ and $R'$ in $\Delta^{\gamma}_{0}$ with $R' \subset Q_{0}(R)^{c}$. Fix the
points $x \in R$ and $y \in R'$. Introduce the shorthand notation
$t_{m} := t_{m}(x,y)$. Define
$$
M := 8 \cdot t_{m} \sqrt{\frac{2^{j(R) + j(R')}}{\abs{x}^{2} + \abs{y}^{2}}}.
$$
 Then we must have the bound
\begin{equation}
\label{eq:FarKernelResult}
k_{t}(x,y) \cdot k_{t_{m}}(x,y)^{-1} \lesssim \frac{1}{2^{(j(R) +
    j(R'))(d+1)}}
\end{equation}
for all $t \leq t_{m} / M = \frac{1}{8} \sqrt{\frac{\abs{x}^{2} +
    \abs{y}^{2}}{2^{j(R) + j(R')}}}$.
 \end{lem}

\begin{proof}  
According to Lemma \ref{lem:Taylor}, $t_{m} / M \leq t_{m}$. 
As $t \mapsto k_{t}(x,y)$ is increasing for $t \leq t_{m}(x,y)$, it
follows that it is sufficient to show \eqref{eq:FarKernelResult} for
the value $t_{m} / M$.
We then have
$$
k_{t_{m}/M}(x,y) \cdot k_{t_{m}}(x,y)^{-1} = M^{d/2} \exp
\br{(\alpha(t_{m}) - \alpha(t_{m}/M))\br{\abs{x}^{2} +
    \abs{y}^{2}}} \cdot \exp \br{- \frac{\abs{x - y}^{2}}{2 t_{m}}(M - 1)}.
$$
Let's find a bound on the function $\alpha(t_{m}) - \alpha(t_{m}/M)$ in
terms of $t_{m}$ and $M$. Define the function $\beta : \R_{>0} \rightarrow \R$ through
$$
\beta(u) := \alpha\br{\frac{1}{u}} = \frac{\sqrt{1 + \frac{1}{u^{2}}}
  - 1}{2/u} = \frac{\sqrt{1 + u^{2}} - u}{2}.
$$
For any $u \leq 1$, perform a Taylor expansion about the origin
for $\beta$ to obtain
$$
\beta(u) = \frac{1}{2} \br{1 - u + \frac{u^{2}}{2} - \frac{u^{4}}{8} +
\frac{u^{6}}{16} - \cdots}.
$$
According to Lemma \ref{lem:Taylor}, both $t_{m}$ and $t_{m}/M$ are
greater than 1. The above formula will therefore apply to 
these values.
$$
\alpha(t_{m}) = \beta(1/t_{m}) = \frac{1}{2} \br{1 - \frac{1}{t_{m}} +
\frac{1}{2 t_{m}^{2}} - \frac{1}{8 t_{m}^{4}} + \frac{1}{16 t_{m}^{6}}
- \cdots},
$$
$$
\alpha(t_{m}/M) = \frac{1}{2} \br{1 - \frac{M}{t_{m}} + \frac{M^{2}}{2
  t_{m}^{2}} - \frac{M^{4}}{8 t_{m}^{4}} + \frac{M^{6}}{16 t_{m}^{6}}
- \cdots}.
$$
Which gives
\begin{align*}\begin{split}   
\alpha(t_{m}) - \alpha(t_{m}/M) &= \frac{(M-1)}{2 t_{m}} - \frac{(M^{2}
  - 1)}{4 t_{m}^{2}} + \frac{(M^{4} - 1)}{16 t_{m}^{4}} - \frac{(M^{6}
  - 1)}{32 t_{m}^{6}} + \cdots \\
&\leq \frac{(M - 1)}{2 t_{m}} - \frac{(M^{2} - 1)}{4 t_{m}^{2}} +
\frac{(M^{4} - 1)}{ 16 t_{m}^{4}}.
 \end{split}\end{align*}
As $M^{2} - 1 \geq \frac{M^{2}}{2}$ and $\frac{(M^{4} - 1)}{16
  t_{m}^{4}} \leq \frac{M^{2}}{16 t_{m}^{2}}$, we obtain
\begin{equation}
\label{eq:FarKernel1}
\alpha(t_{m}) - \alpha(t_{m} / M) \leq \frac{(M - 1)}{2 t_{m}} -
\frac{M^{2}}{16 t_{m}^{2}}.
\end{equation}
Once more from Lemma \ref{lem:Taylor}, we have that
\begin{align*}\begin{split}  
 M^{d/2} 2^{(j(R) + j(R'))(d+1)} &\leq t_{m}^{d/2} 2^{(j(R) +
   j(R'))(d+1)} \\
&\lesssim \abs{y - x}^{d/2}
 2^{(j(R) + j(R'))(d+1)} \\
&\leq (\abs{y} + \abs{x})^{d/2} 2^{(j(R) + j(R'))(d+1)} \\
&\lesssim (2^{j(R)} + 2^{j(R')})^{d/2} 2^{(j(R)+j(R'))(d+1)}.
 \end{split}\end{align*}

It is easy to see that there must exist some $A \geq 0$, independent
of both $R$ and $R'$, such that
$$
\br{2^{j(R)} + 2^{j(R')}}^{d/2} 2^{(j(R) + j(R'))(d+1)} \leq A
e^{2^{j(R) + j(R')}}.
$$
This would then give
$$
M^{d/2} 2^{(j(R) + j(R'))(d+1)} \lesssim  e^{2^{j(R) + j(R')}}.
$$
On applying \eqref{eq:FarKernel1} and the above,
\begin{align*}\begin{split}  
&k_{t_{m}/M}(x,y) \cdot k_{t_{m}}(x,y)^{-1} 2^{(j(R) + j(R'))(d+1)}
 \\
& \ \ \lesssim  M^{d/2} 2^{(j(R) + j(R'))(d+1)} \exp \br{(\alpha(t_{m}) -
   \alpha(t_{m}/M)) \br{\abs{x}^{2} + \abs{y}^{2}}} \cdot \exp
 \br{-\frac{\abs{x-y}^{2}}{2 t_{m}} (M - 1)} \\ & \ \ \lesssim \exp
 \br{2^{j(R) + j(R')}} \cdot \exp \br{\br{\frac{(M-1)}{2 t_{m}} -
     \frac{M^{2}}{16 t_{m}^{2}}}\br{\abs{x}^{2} + \abs{y}^{2}}} \cdot
 \exp \br{- \frac{\abs{x - y}^{2}}{2 t_{m}}(M-1)} \\
& \ \ = \exp \br{2^{j(R) + j(R')} + \frac{(M-1)}{t_{m}} \langle x , y \rangle -
  \frac{M^{2}}{16 t_{m}^{2}} \br{\abs{x}^{2} + \abs{y}^{2}}} \\
& \ \ \leq \exp \br{2^{j(R) + j(R')} + \frac{(M-1)}{t_{m}} \abs{x} \abs{y} -
  \frac{M^{2}}{16 t_{m}^{2}} \br{\abs{x}^{2} + \abs{y}^{2}}}.
 \end{split}\end{align*}
On applying $M / t_{m} \leq 1 / \br{2^{4} d^{2}}$,
\begin{align*}\begin{split}  
 k_{t_{m} /M}(x,y) \cdot k_{t_{m}}(x,y)^{-1} & 2^{(j(R) + j(R'))(d+1)} \\
&\lesssim \exp \br{2^{j(R) + j(R')} + \frac{\abs{x} \abs{y}}{2^{4}
    d^{2}} - \frac{M^{2}}{16 t_{m}^{2}} \br{\abs{x}^{2} +
    \abs{y}^{2}}} \\
&\lesssim \exp \br{2^{j(R) + j(R')} + \frac{2^{j(R) + j(R')}}{2^{4} d}
  - \frac{M^{2}}{16 t_{m}^{2}} \br{\abs{x}^{2} + \abs{y}^{2}}}.
 \end{split}\end{align*}
From which the definition of $M$ then provides
$$
k_{t_{m}/M}(x,y) \cdot k_{t_{m}}(x,y)^{-1} 2^{(j(R) + j(R'))(d+1)}
\lesssim 1.
$$
 \end{proof}

The next result is a direct analogue for $A_{p}^{+}$ of the defining
condition for the classic $A_{p}$ class. It is unlikely that this
condition is enough to completely characterise $A_{p}^{+}$.

\begin{lem} 
 \label{lem:ApProperty}
Let $w$ be a weight on $\R^{d}$ and suppose that
$\mathcal{M}^{+}_{far} : L^{p}(w) \rightarrow L^{p}(w)$ is
bounded for some $1 < p < \infty$. Fix cubes $R$ and $R'$ in $\Delta^{\gamma}_{0}$ with $R' \not\subset
Q_{0}(R)$. Then there must exist some constant $C > 0$, independent of
both $R$ and $R'$, such that
\begin{equation}
\label{eq:lem:ApProperty}
w(R)^{\frac{1}{p}} \cdot w^{-\frac{1}{p-1}}(R')^{\frac{p-1}{p}} \leq C
\cdot k_{t_{m}(\tilde{x},\tilde{y})}(\tilde{x},\tilde{y})^{-1}
\end{equation}  
for all $\tilde{x} \in R$ and $\tilde{y} \in R'$.
 \end{lem}

\begin{proof}  
 It shall first be shown that
$$
R' \subset Q_{t_{m}(\tilde{x}, \tilde{y})}(R).
$$
Fix any point $y \in R'$. From the definition of $Q_{t}(R)$, it will
be sufficient to show that
$$
\abs{y} \leq  2^{2} t_{m}(\tilde{x}, \tilde{y})^{2} 2^{j(R)}. 
$$
First suppose that $R$ is not in the first layer. Then
\begin{align*}\begin{split}  
 2^{j(R)} \cdot t_{m}(\tilde{x}, \tilde{y})^{2} &\geq  2^{j(R)}
 \frac{\abs{\tilde{y}}^{2}}{9^{2} d^{2} \abs{\tilde{x}}^{2}}.
 \end{split}\end{align*}
As $\abs{\tilde{x}} \leq \sqrt{d} \cdot 2^{j(R)}$, $\abs{y} \leq
\sqrt{d} \cdot 2
\abs{\tilde{y}}$ and $\abs{\tilde{y}} \geq d^{4} 2^{16} 2^{j(R)}$, we
have that
\begin{align*}\begin{split}  
 2^{j(R)} \cdot t_{m}(\tilde{x}, \tilde{y})^{2} &\geq
\frac{2^{j(R)}}{9^{2} d^{2}} \cdot \frac{\abs{y}}{2 \sqrt{d}} \cdot \frac{d^{4} 2^{16}
  2^{j(R)}}{d 2^{2j(R)}} \\
&\geq \abs{y}.
 \end{split}\end{align*}
Next suppose that $R$ is contained in the first layer. Then
\begin{align*}\begin{split}  
 2^{j(R)} \cdot t_{m}(\tilde{x}, \tilde{y})^{2} &\geq
 \frac{\abs{\tilde{y}}^{2}}{9^{2} d^{2}} \\
&\geq \frac{\abs{y}}{2 \sqrt{d}} \cdot \frac{2^{16} d^{4}}{9^{2} d^{2}} \\
&\geq \abs{y}. 
 \end{split}\end{align*}
This demonstrates that $R' \subset Q_{t_{m}(\tilde{x},
  \tilde{y})}(R)$. Then, for any $\tilde{x} \in R$ and $\tilde{y} \in R'$,
\begin{align*}\begin{split}  
 w(R) \br{\int_{R'} \abs{f(y)} dy}^{p} &= \int_{R} w(x) dx
 \frac{k_{t_{m}(\tilde{x}, \tilde{y})}(\tilde{x},
   \tilde{y})^{p}}{k_{t_{m}(\tilde{x}, \tilde{y})}(\tilde{x},
   \tilde{y})^{p}} \br{ \int_{R'} \abs{f(y)} dy}^{p} \\
&=\frac{1}{k_{t_{m}(\tilde{x}, \tilde{y})}(\tilde{x}, \tilde{y})^{p}}
\int_{R} \br{k_{t_{m}(\tilde{x}, \tilde{y})}(\tilde{x}, \tilde{y})
  \int_{R'} \abs{f(y)} dy}^{p} w(x) dx \\
&\leq \frac{1}{k_{t_{m}(\tilde{x}, \tilde{y})}(\tilde{x},
  \tilde{y})^{p}} \int_{R} \mathcal{M}^{+}_{far}(f \cdot \chi_{R'})(x)^{p} w(x) dx.
 \end{split}\end{align*}
From the boundedness of $\mathcal{M}^{+}_{far}$, we then obtain
$$
w(R) \br{\int_{R'} \abs{f(y)} dy}^{p} \lesssim
\frac{1}{k_{t_{m}(\tilde{x}, \tilde{y})}(\tilde{x}, \tilde{y})^{p}}
\int_{R'} \abs{f(y)}^{p} w(y) dy.
$$
Take $f := \br{w + \varepsilon}^{-\frac{1}{p-1}}$ for some
$\varepsilon > 0$. Then
$$
w(R) \br{ \int_{R'} \br{w(y) + \varepsilon}^{-\frac{1}{p-1}} dy}^{p}
\lesssim \frac{1}{k_{t_{m}(\tilde{x}, \tilde{y})}(\tilde{x},
  \tilde{y})^{p}} \int_{R'} \frac{w(y)}{\br{w(y) +
    \varepsilon}^{\frac{p}{p-1}}} dy
$$
for all $\varepsilon > 0$.  Which implies that
$$
w(R) \br{ \int_{R'} \br{w(y) + \varepsilon}^{-\frac{1}{p-1}} dy}^{p}
\lesssim \frac{1}{k_{t_{m}(\tilde{x}, \tilde{y})}(\tilde{x},
  \tilde{y})^{p}} \int_{R'} \frac{\br{w(y) + \varepsilon}}{\br{w(y) +
    \varepsilon}^{\frac{p}{p-1}}} dy
$$
$$
\Rightarrow \quad w(R) \br{\int_{R'} \br{w(y) +
    \varepsilon}^{-\frac{1}{p-1}} dy}^{p-1} \lesssim
\frac{1}{k_{t_{m}(\tilde{x}, \tilde{y})}(\tilde{x}, \tilde{y})^{p}}
$$
for each $\varepsilon > 0$. An application of the Lebesgue monotone
convergence theorem then produces the desired result.
\phantom{This is empty space}
\end{proof}

Finally, enough machinery is in place to prove our main result.

\begin{customthm}{A}
 \label{thm:far1}
Let $w$ be a weight on $\R^{d}$ and $1 < p < \infty$. Then we have
\begin{equation}
\label{eq:prop:far1}
\norm{\mathcal{M}_{far}^{+}}_{L^{p}(w)} < \infty 
\quad \Rightarrow \quad \norm{\mathcal{T}^{*}_{far}}_{L^{p}(w)} < \infty \quad \Rightarrow \quad
\norm{\mathcal{M}^{-}_{far}}_{L^{p}(w)} < \infty.
\end{equation} 
 \end{customthm}

\begin{proof}  
The second implication follows quickly from the pointwise
 bound
\begin{align*}\begin{split}  
 \mathcal{T}^{*}_{far}f(x) &= \sup_{t > 0} \int_{F(R)} k_{t}(x,y)
 \abs{f(y)} dy \\
&= \sup_{t > 0} \sum_{R' \in \mathcal{F}(R)} \int_{R'} k_{t}(x,y)
\abs{f(y)} dy \\
&\geq \sup_{t > 0} \sum_{R' \in \mathcal{F}(R), \ R' \subset Q_{t}(R)}
k_{t}^{-}(R,R')
\int_{R'} \abs{f(y)} dy \\
&= \mathcal{M}^{-}_{far}f(x)
 \end{split}\end{align*}
for any $f \in L^{1}_{loc}(\R^{d})$ and $x \in R \in
\Delta^{\gamma}_{0}$.

As for the first implication, suppose that $\norm{\mathcal{M}^{+}_{far}}_{L^{p}(w) \rightarrow
  L^{p}(w)} < \infty$. Then
\begin{align*}\begin{split}  
 \norm{\mathcal{T}^{*}_{far}f}_{L^{p}(w)} &= \brs{ \int_{\R^{d}}
 \abs{\mathcal{T}^{*}_{far}f(x)}^{p} w(x) dx }^{1/p}  \\
 &= \brs{ \int_{\R^{d}} \br{ \sup_{t > 0}
   e^{- t \mathcal{L}} \br{f \cdot \chi_{F(R_{x})}}(x)}^{p} w(x) dx. }^{1/p} \\
&= \brs{ \int_{\R^{d}} \br{ \sup_{t > 0} \int_{F(R_{x})} k_{t}(x,y) \abs{f(y)} dy }^{p} w(x) dx }^{1/p}.
 \end{split}\end{align*}
The heat operators can be expanded dyadically to obtain
\begin{align*}\begin{split}  
 \norm{\mathcal{T}^{*}_{far}f}_{L^{p}(w)} &= \brs{ \int_{\R^{d}} \br{
     \sup_{t > 0} \sum_{R' \in \mathcal{F}(R_{x})} \int_{R'} k_{t}(x,y) \abs{f(y)} dy}^{p} w(x) dx}^{1/p} \\
&\lesssim  \brs{ \int_{\R^{d}} \br{
     \sup_{t > 0} \sum_{R' \in \mathcal{F}(R_{x})} k_{t}^{+}(R_{x}, R') \norm{f}_{L^{1}(R')}}^{p} w(x) dx}^{1/p} \\
&\lesssim  \brs{ \int_{\R^{d}} \br{
     \sup_{t > 0} \sum_{R' \in \mathcal{F}(R_{x}), \ R' \subset Q_{t}(R_{x})
       } k_{t}^{+}(R_{x}, R') \norm{f}_{L^{1}(R')} \right. \right. \\
& \qquad \qquad \left. \left. +  \sup_{t > 0} \sum_{R' \in \mathcal{F}(R_{x}), \ R' \not\subset Q_{t}(R_{x})
       } k_{t}^{+}(R_{x}, R') \norm{f}_{L^{1}(R')}}^{p} w(x) dx}^{1/p}.
 \end{split}\end{align*}
On applying Minkowski's inequality,
\begin{align*}\begin{split}  
 \norm{\mathcal{T}^{*}_{far}f}_{L^{p}(w)} &\lesssim \brs{\int_{\R^{d}}
 \br{\sup_{t > 0} \sum_{R' \in \mathcal{F}(R_{x}), \ R' \subset Q_{t}(R_{x})} k_{t}^{+}(R_{x}, R') \norm{f}_{L^{1}(R')}}^{p}
 w(x) dx}^{1/p} \\ & \qquad + \brs{\int_{\R^{d}}
 \br{\sup_{t > 0} \sum_{R' \in \mathcal{F}(R_{x}), \ R' \not\subset Q_{t}(R_{x})} k_{t}^{+}(R_{x}, R') \norm{f}_{L^{1}(R')}}^{p}
 w(x) dx}^{1/p} \\
&= \norm{\mathcal{M}^{+}_{far} f}_{L^{p}(w)} + \brs{\int_{\R^{d}}
 \br{\sup_{t > 0} \sum_{R' \in \mathcal{F}(R_{x}), \ R' \not\subset Q_{t}(R_{x})} k_{t}^{+}(R_{x}, R') \norm{f}_{L^{1}(R')}}^{p} w(x) dx}^{1/p}.
 \end{split}\end{align*}
It remains to bound the tail end term on the right hand side of the
above expression. On expanding dyadically once more,
\begin{align*}\begin{split}  
 \int_{\R^{d}}
 &\br{\sup_{t > 0} \sum_{R' \in \mathcal{F}(R_{x}), \ R' \not\subset Q_{t}(R_{x})} k_{t}^{+}(R_{x}, R') \norm{f}_{L^{1}(R')}}^{p}
 w(x) dx \\   & \qquad \qquad = \sum_{R \in \Delta^{\gamma}_{0}} \int_{R}  \br{\sup_{t > 0} \sum_{R' \in \mathcal{F}(R), \ R' \not\subset Q_{t}(R)} k_{t}^{+}(R, R') \norm{f}_{L^{1}(R')}}^{p}
 w(x) dx \\
 & \qquad \qquad = \sum_{R \in \Delta^{\gamma}_{0}}  \br{\sup_{t > 0} \sum_{R' \in \mathcal{F}(R), \ R' \not\subset Q_{t}(R)} k_{t}^{+}(R, R')
   \norm{f}_{L^{1}(R')} w(R)^{1/p}}^{p}.
 \end{split}\end{align*}
Let $x_{R}^{t}$ and $y_{R'}^{t}$ denote points contained in
$R$ and $R'$ respectively that satisfy
$$
k_{t}^{+}(R,R') \leq 2 \cdot k_{t}(x_{R}^{t}, y_{R'}^{t}).
$$
On applying H\"{o}lder's property and Lemma \ref{lem:ApProperty} we obtain
\begin{align*}\begin{split}  
 \sum_{R \in \Delta^{\gamma}_{0}} & \br{\sup_{t > 0} \sum_{R'
     \in \mathcal{F}(R), \ R' \not\subset Q_{t}(R)}
   k_{t}^{+}(R,R') \norm{f}_{L^{1}(R')} w(R)^{1/p}}^{p} \\ & \qquad \qquad \lesssim
\sum_{R \in \Delta^{\gamma}_{0}} \br{\sup_{t > 0} \sum_{R'
     \in \mathcal{F}(R), \ R' \not\subset Q_{t}(R)}
   k_{t}^{+}(R,R') w^{-\frac{1}{p-1}}(R')^{\frac{p-1}{p}}
   w(R)^{\frac{1}{p}} \norm{f}_{L^{p}(R',w)}}^{p} \\
& \qquad \qquad \lesssim 
\sum_{R \in \Delta^{\gamma}_{0}} \br{\sup_{t > 0} \sum_{R'
     \in \mathcal{F}(R), \ R' \not\subset Q_{t}(R)}
   k_{t}(x_{R}^{t},y_{R'}^{t}) \cdot k_{t_{m}(x_{R}^{t},y_{R}^{t})}(x_{R}^{t},y_{R'}^{t})^{-1} \norm{f}_{L^{p}(R',w)}}^{p}.
 \end{split}\end{align*}
Note that since $\abs{y^{t}_{R'}} \geq 2^{8} t^{2} 2^{j(R)}$, it
follows that
\begin{align*}\begin{split}  
\frac{1}{8} \sqrt{\frac{\abs{x_{R}^{t}}^{2} + \abs{y_{R'}^{t}}^{2}}{2^{j(R) + j(R')}}}
&\geq \frac{1}{8} \sqrt{\frac{\abs{y_{R'}^{t}}^{2}}{2^{j(R) + j(R')}}} \\
&\geq \frac{1}{8} \sqrt{\frac{2^{j(R') - 1} \cdot 2^{8} t^{2} 2^{j(R)}}{2^{j(R) +
      j(R')}}} \\
&\geq t. 
 \end{split}\end{align*}
This implies that Lemma \ref{lem:FarKernel} can be applied to obtain
\begin{align*}\begin{split}  
 \sum_{R \in \Delta^{\gamma}_{0}}& \br{\sup_{t > 0} \sum_{R' \in
  \mathcal{F}(R), \ R' \not\subset Q_{t}(R)}
k_{t}(x_{R}^{t}, y_{R'}^{t}) \cdot k_{t_{m}(x_{R}^{t}, y_{R'}^{t})}(x_{R}^{t},
y_{R'}^{t})^{-1} \norm{f}_{L^{p}(R',w)}}^{p}\\ & \qquad \qquad \lesssim  \sum_{R \in \Delta^{\gamma}_{0}} \br{\sup_{t > 0} \sum_{R' \in
  \mathcal{F}(R), \ R' \not\subset Q_{t}(R)}
2^{-(j(R) + j(R'))(d+1)} \norm{f}_{L^{p}(R',w)}}^{p} \\
& \qquad \qquad \lesssim \norm{f}_{L^{p}(w)}^{p} \sum_{k=0}^{\infty} \sum_{R \in
  L_{k}} \br{\sum_{l=0}^{\infty} \sum_{R' \in L_{l}}
  2^{-(k+l)(d+1)}}^{p} \\
& \qquad \qquad \lesssim \norm{f}_{L^{p}(w)}^{p} \sum_{k = 0}^{\infty} 2^{kd}
\br{\sum_{l = 0}^{\infty} 2^{ld} \cdot 2^{-(k+l)(d+1)}}^{p} \\
& \qquad \qquad \lesssim \norm{f}_{L^{p}(w)}^{p},
 \end{split}\end{align*}
since the number of cubes in a layer $L_{k}$ is bounded by a constant
multiple of $2^{kd}$.
 \end{proof}

Theorems \ref{thm:Main} and
\ref{thm:Main2}, together with the fact that $\norm{\mathcal{T}^{*}}_{L^{p}(w)} <
\infty$ if and only if both $\norm{\mathcal{T}^{*}_{loc}}_{L^{p}(w)
  \rightarrow L^{p}(w)} < \infty$ and
$\norm{\mathcal{T}^{*}_{far}}_{L^{p}(w) \rightarrow L^{p}(w)} < \infty$ for
any weight $w$ on $\R^{d}$, lead to the below corollary.

\begin{cor}
\label{cor:Classes}
The following chain of inclusions holds for any $1 < p < \infty$,
\begin{equation}
\label{eq:Classes}
A_{p}^{+} \subseteq \lb w \ \mathrm{weight} \ \mathrm{on} \ \R^{d} : \norm{\mathcal{T}^{*}}_{L^{p}(w)
  \rightarrow L^{p}(w) < \infty} \rb \subseteq A_{p}^{-}.
\end{equation}
\end{cor}

The class of weights in the middle of the above chain of inclusions is
a natural candidate for the $A_{p}$ class associated with the harmonic
oscillator. The above corollary indicates that our $A_{p}$ classes are
honing in on what should be the correct class.

\section{Relation to the $A_{p}^{\infty}$ Class}
\label{sec:final}

Recall the definitions of the classes $A_{p}^{\infty}$ and
$A_{p}^{\theta}$ from section \ref{sec:intro}. This section is devoted
to the proof of the strict inclusion $A_{p}^{\infty} \subsetneq A_{p}^{+}$. This
will be accomplished by first showing, for any $\theta \geq 0$, that
the
pointwise bound $\mathcal{M}^{+}_{far}f(x) \lesssim M^{\theta}f(x)$
holds for
all $f \in L^{1}_{loc}(\R^{d})$ and $x \in \R^{d}$, thereby
demonstrating the inclusion $A_{p}^{\theta} \subseteq
A_{p}^{far+}$. The following upper bound for the heat kernel $k$ will
be utilised. Refer to \cite{kurata2000estimate} for proof. 

\begin{lem} 
 \label{lem:bound}
For any $N > 0$, there exists a constant $C_{N} > 0$ such that
 \begin{equation}
\label{eq:bound}
k_{\sinh 2 t}(x,y) \leq C_{N} t^{-d/2} \exp \br{-\frac{\abs{x-y}^{2}}{2 t}}
\br{1 + \frac{\sqrt{t}}{\rho(x)} + \frac{\sqrt{t}}{\rho(y)}}^{-N}
\end{equation}
for all $x, y \in \R^{d}$.
 \end{lem}

Recall that the $\sinh 2 t$ factor in the above expression is due to
the kernel
rescaling introduced in section \ref{sec:local}.

\begin{prop} 
 \label{prop:Pointwise} 
 For any $\theta \geq 0$, there exists some $C_{\theta} > 0$ so
 that
$$
\mathcal{M}^{+}_{far}f(x) \leq C_{\theta} M^{\theta} f(x)
$$ 
for every locally integrable function $f$ on $\R^{d}$ and $x \in
\R^{d}$.
 \end{prop}

\begin{proof}  
For $R \in \Delta^{\gamma}_{0}$ and $k \geq 0$, define
$\mathcal{C}_{k}(R)$ to be the collection of cubes
$R' \in \Delta^{\gamma}_{0}$ that satisfy $d(R,R') < 2^{k} l(R)$.
 As $\mathcal{F}(R)
\subset \Delta^{\gamma}_{0} / \mathcal{C}_{0}(R)$, the operator
$\mathcal{M}^{+}_{far}$ can be decomposed as
\begin{align*}\begin{split}  
 \mathcal{M}^{+}_{far}f(x) &\leq  \sup_{t > 0} \sum_{R' \in
 \Delta^{\gamma}_{0} / \mathcal{C}_{0}(R)}
 k_{t}^{+}(R,R') \int_{R'} \abs{f(y)} dy \\
&= \sup_{t > 0} \sum_{R' \in
 \Delta^{\gamma}_{0} / \mathcal{C}_{0}(R)}
 k_{ \sinh 2 t}^{+}(R,R') \int_{R'} \abs{f(y)} dy \\
&\leq \sup_{t > 0} \sum_{k = 1}^{\infty} \sum_{R' \in
  \mathcal{C}_{k}(R) / \mathcal{C}_{k-1}(R)} k_{\sinh 2 t}^{+}(R,R') \int_{R'}
\abs{f(y)} dy
 \end{split}\end{align*}
for $x \in R$. Let's find a bound on the values $k_{\sinh 2 t}^{+}(R,R')$ for $R'
\subset \mathcal{C}_{k}(R) / \mathcal{C}_{k-1}(R)$.  Suppose that $x
\in R$ and $y \in R' \in \mathcal{C}_{k}(R) / \mathcal{C}_{k-1}(R)$
where $k \geq 1$. Then, $\abs{x - y} \geq 2^{k-1} 2^{-j(R)}$.
From this bound, Lemma \ref{lem:bound} and the inequality $\rho(x)
\leq 2^{1 - j(R)}$,
\begin{align*}\begin{split}  
 k_{\sinh 2 t}(x,y) &\lesssim t^{-d/2} \exp \br{-\frac{\abs{x-y}^{2}}{2t}}
 \br{1 + \frac{\sqrt{t}}{\rho(x)}}^{-N} \\
&\lesssim t^{-d/2} \frac{t^{M/2}}{\abs{x - y}^{M}} \br{1 +
  2^{j(R) - 1} \sqrt{t}}^{-N} \\
&\lesssim t^{-d/2} \br{2^{j(R)} \sqrt{t}}^{M} 2^{-k M} \br{1 +
  2^{j(R) - 1} \sqrt{t}}^{-N} \\
&\lesssim 2^{j(R) d} 2^{-k M} \br{2^{j(R) - 1} \sqrt{t}}^{M-d} \br{1 + 2^{j(R) - 1} \sqrt{t}}^{-N}
 \end{split}\end{align*}
for any $M > 0$. Therefore
\begin{equation}
\label{eq:Inclusion1}
k_{\sinh 2 t}^{+}(R,R') \lesssim 2^{j(R) d} 2^{-k M} \br{
    2^{j(R) - 1} \sqrt{t}}^{M-d} \br{1 + 2^{j(R) - 1} \sqrt{t}}^{-N}
\end{equation}
for any $R' \subset \mathcal{C}_{k}(R) / \mathcal{C}_{k-1}(R)$. On applying this bound to our
previous decomposition we find that $\mathcal{M}^{+}_{far}f(x)$ can be
estimated above by
$$
\sup_{t > 0 } \sum_{k=1}^{\infty} 2^{j(R)
   d} 2^{-k M} \br{2^{j(R) - 1} \sqrt{t}}^{M-d} \br{1 +
   2^{j(R) - 1} \sqrt{t}}^{-N} \sum_{R' \in
   \mathcal{C}_{k}(R) / \mathcal{C}_{k-1}(R)} \int_{R'} \abs{f(y)} dy.
$$
Define $R_{k}$ to be the smallest cube that contains
every cube in the collection $\mathcal{C}_{k}(R)$. Then
$$
\mathcal{M}^{+}_{far}f(x) \lesssim 2^{j(R) d} \sup_{s > 0} s^{M-d} \br{1 + s}^{-N}
\sum_{k=1}^{\infty} 2^{-kM} \int_{R_{k}} \abs{f(y)} dy,
$$
where we set $s := 2^{j(R) - 1} \sqrt{t}$. It is obvious
that if we set $N \geq M - d$, then the supremum term must be bounded
by $1$. We then  obtain
$$
\mathcal{M}^{+}_{far}f(x) \lesssim 2^{j(R) d} \sum_{k=1}^{\infty} 2^{-k M} \int_{R_{k}}
\abs{f(y)} dy.
$$
On noting that $l(R_{k}) \approx 2^{k} 2^{-j(R)}$ and
$\psi_{\theta}(R_{k}) \lesssim 2^{k \theta}$, we have
\begin{align*}\begin{split}  
 \mathcal{M}^{+}_{far}f(x) &\lesssim \sum_{k=1}^{\infty} 2^{-k(M - d - \theta)}
 \frac{1}{2^{k \theta}} \br{\frac{2^{j(R)d}}{2^{k d}}} \int_{R_{k}}
 \abs{f(y)} dy \\
&\lesssim \sum_{k=1}^{\infty} 2^{-k(M - d - \theta)}
\frac{1}{\psi_{\theta}(R_{k}) \abs{R_{k}}}
\int_{R_{k}}\abs{f(y)} dy \\
&\leq M_{\theta}f(x)
 \end{split}\end{align*}
for $M \geq d + \theta$.
 \end{proof}

\begin{customprop}{C}
 \label{prop:StrictInclusion} 
 The following chain of strict inclusions holds for any $1 < p < \infty$,
$$
A_{p} \subsetneq A_{p}^{\infty} \subsetneq A_{p}^{+}.
$$
 \end{customprop}

\begin{proof}  
The strict inclusion $A_{p} \subsetneq A_{p}^{\infty}$ has already
been proved in \cite{bongioanni2011classes}. As for the upper inclusion, the previous proposition demonstrates that $A_{p}^{\infty} \subseteq
 A_{p}^{far+}$. It will now be proved that $A_{p}^{\infty}
 \subseteq A_{p}^{loc}$. Fix $w \in A_{p}^{\infty}$. Then there must
 exist some $\theta \geq 0$ such that $w \in A_{p}^{\theta}$. It must
 be shown that there exists some $B > 0$ that satisfies
\begin{equation}
\label{eq:uniformlocal}
\brs{w}_{A_{p}(N(R))} \leq B
\end{equation}
for every $R \in \Delta^{\gamma}_{0}$. Fix any cube $R \in
\Delta^{\gamma}_{0}$ and $Q$ a dyadic subcube of $N(R)$. As $w \in A_{p}^{\theta}$,
there must exist some $C > 0$ such that
$$
w(Q)^{\frac{1}{p}} w^{-\frac{1}{p-1}}(Q)^{\frac{p-1}{p}} \leq C
\abs{Q} \br{1 + \frac{l(Q)}{\rho(c_{Q})}}^{\theta}.
$$
As $Q$ is a dyadic subcube of $N(R)$, we have that $l(Q) \leq 4
\rho(c_{R})$ and $\rho(c_{Q}) \geq \rho(c_{R}) / 2$. Therefore
\begin{align*}\begin{split}  
 w(Q)^{\frac{1}{p}} w^{-\frac{1}{p-1}}(Q)^{\frac{p-1}{p}} &\leq C
 \abs{Q} \br{1 + 8}^{\theta} \\
&\leq 9^{\theta} C \abs{Q}.
 \end{split}\end{align*}
This demonstrates that \eqref{eq:uniformlocal} holds with constant $B
:= 9^{\theta} C$.

It will now be proved that the inclusion of $A_{p}^{\infty}$ in
$A_{p}^{+}$ is in fact strict. In particular, the weight defined by
$$
w(x) = w(x_{1}, \cdots, x_{d}) = e^{\abs{x_{1}}}
$$
for $x \in \R^{d}$ will be shown to belong to the class $A_{p}^{+}$
but not $A_{p}^{\infty}$.

Let's first show that $w \in A_{p}^{loc}$. That is, it will be proved that there exists $C > 0$ such that for any $R \in \Delta^{\gamma}_{0}$ and dyadic
subcube $Q$ of $N(R)$
\begin{equation}
  \label{eqtn:LocalExp}
w(Q) w^{-\frac{1}{p-1}}(Q)^{p-1} \leq C \abs{Q}^{p}.
\end{equation}
Note that for any $x = (x_{1}, \cdots, x_{d}) \in Q$ we must have the bound
$$
\abs{c_{R}^{(1)}} - 4 \cdot 2^{-j(R)} \leq \abs{x_{1}} \leq
\abs{c^{(1)}_{R}} + 4 \cdot 2^{-j(R)},
$$
where $c_{R} = \br{c_{R}^{(1)}, \cdots, c_{R}^{(d)}}$. This gives
\begin{align*}\begin{split}  
    w(Q) &= \int_{Q} e^{\abs{x_{1}}} dx \\
    &\lesssim  e^{\abs{c_{R}^{(1)}}} \abs{Q}.
 \end{split}\end{align*}
Similarly,
\begin{align*}\begin{split}  
 w^{-\frac{1}{p-1}}(Q)^{p-1} &= \br{\int_{Q} e^{-
     \frac{\abs{x_{1}}}{p-1}} d x}^{p-1} \\
 &\lesssim e^{- \abs{c_{R}^{(1)}}} \abs{Q}^{p-1}.
\end{split}\end{align*}
This gives estimate \eqref{eqtn:LocalExp} and proves that $w \in A_{p}^{loc}$.

Next let's prove that $w \in A_{p}^{far+}$. That is, it must be shown
that
$$
\norm{\mathcal{M}^{+}_{far}f}_{L^{p}(w)} \lesssim \norm{f}_{L^{p}(w)}
$$
for any $f \in L^{p}(w)$.
\begin{align*}\begin{split}  
 \norm{\mathcal{M}^{+}_{far} f}^{p}_{L^{p}(w)} &= \int_{\R^{d}}
 \mathcal{M}^{+}_{far} f(x)^{p} w(x) dx \\
 &= \int_{\R^{d}} \br{\sup_{t > 0} \sum_{\substack{R' \in
       \mathcal{F}(R_{x}), \\ R' \subset Q_{t}(R_{x})}} k_{t}^{+}(R_{x}, R')
   \int_{R'} \abs{f(y)} dy}^{p} w(x) dx \\
 &= \int_{\R^{d}} \sup_{t > 0} \br{\int_{Q_{t}(R_{x}) \cap
     F(R_{x})} k_{t}^{+}(R_{x}, R_{y}) \abs{f(y)}
 w(y)^{\frac{1}{p}} w(y)^{- \frac{1}{p}} dy}^{p} w(x) dx.
\end{split}\end{align*}
On applying H\"{o}lder's inequality we obtain
\begin{align}\begin{split}  
    \label{eqtn:PropC1}
      \norm{\mathcal{M}^{+}_{far} f}^{p}_{L^{p}(w)} &\lesssim \br{
  \int_{\R^{d}} \sup_{t > 0} \br{\int_{Q_{t}(R_{x}) \cap
      F(R_{x})} k_{t}^{+}(R_{x}, R_{y})^{p'} w(y)^{-\frac{p'}{p}} dy}^{\frac{p}{p'}}
  w(x) dx} \norm{f}^{p}_{L^{p}(w)} \\
&\leq    \br{
  \int_{\R^{d}} \sup_{t > 0} \br{\int_{Q_{t}(R_{x}) \cap
      F(R_{x})} k_{t}^{+}(R_{x}, R_{y})^{p'} dy}^{\frac{p}{p'}}
  w(x) dx} \norm{f}^{p}_{L^{p}(w)}
 \end{split}\end{align}
Let $M \geq 1$, the exact value to be determined at a later time. It
will now be proved that the function
\begin{equation}
  \label{eqtn:PropC11}
 (t,x) \ \mapsto \ \br{\int_{Q_{t}(R_{x}) \cap F(R_{x})}
  k_{t}^{+}(R_{x},R_{y})^{p'} dy}^{\frac{p}{p'}}
\end{equation}
is uniformly bounded for $t > 0$ and $x \in [-M,M]^{d}$. For $x \in \R^{d}$ and $y \in
\R^{d}$, let $\tilde{x}$ and $\tilde{y}$ denote points in $R_{x}$ and
$R_{y}$ respectively that satisfy $k_{t}(R_{x},R_{y}) \leq 2
k_{t}(\tilde{x}, \tilde{y})$. As $\tilde{y} \in F(R_{x}) =
F(R_{\tilde{x}})$ we must have $\abs{\tilde{x} -\tilde{y}}
\geq 2^{-j(R_{x})}$. This implies that
\begin{align*}\begin{split}  
 k_{t}(R_{x}, R_{y}) &\lesssim \frac{1}{\br{2 \pi t}^{\frac{d}{2}}} \exp
 \br{- \frac{\abs{\tilde{x} -  \tilde{y}}^{2}}{2 t}} \cdot \exp \br{-
   \alpha(t) \br{\abs{\tilde{x}}^{2} + \abs{\tilde{y}}^{2}}} \\
 &\lesssim \frac{1}{t^{\frac{d}{2}}} \exp \br{- \frac{2^{- 2j(R_{x})}}{2
     t}} \\
 &\lesssim \frac{1}{t^{\frac{d}{2}}} \cdot \frac{1}{\br{2^{-2j(R_{x})}/2
     t}^{\frac{d}{2}}} \\
 &\approx 2^{d j(R_{x})}.
\end{split}\end{align*}
As $x$ is restricted to $[-M,M]^{d}$, the layer number $j(R_{x})$ is
bounded implying that $(t,x,y) \mapsto k_{t}(R_{x},R_{y})$ is bounded. For $t \leq 1$ the size of $Q_{t}(R_{x})$ is bounded proving
that \eqref{eqtn:PropC11} is bounded for $t \leq 1$ and $x \in
[-M,M]^{d}$. For $t > 1$ note that
\begin{align*}\begin{split}  
 k_{t}(R_{x},R_{y}) &\lesssim \frac{1}{\br{2 \pi t}^{\frac{d}{2}}}
 \exp \br{- \frac{\abs{\tilde{x} - \tilde{y}}^{2}}{2 t}} \cdot \exp \br{- \alpha (t)
   \br{\abs{\tilde{x}}^{2} + \abs{\tilde{y}}^{2}}} \\
 &\lesssim \exp \br{- \alpha (t) \abs{\tilde{y}}^{2}}.
\end{split}\end{align*}
Since $\abs{y} \leq 2 \br{\abs{\tilde{y}} + \sqrt{d}}$ and $\alpha$
is an increasing function,
$$
 \int_{Q_{t}(R_{x}) \cap F(R_{x})} k_{t}(R_{x}, R_{y})^{p'}
 dy \lesssim  \int_{\R^{d}} \exp \br{- \alpha (1) \br{\frac{\abs{y}}{2} -
     \sqrt{d}}^{2} p'} dy
$$
which is clearly integrable. This shows that \eqref{eqtn:PropC11} is
uniformly bounded for $x \in [-M,M]^{d}$ and $t > 0$. Therefore, to
complete the proof of $w \in A_{p}^{far+}$ it is sufficient to show that
$$
\int_{\R^{d} / [-M,M]^{d}} \sup_{t > 0} \br{\int_{Q_{t}(R_{x}) \cap
    F(R_{x})} k_{t}^{+}(R_{x}, R_{y})^{p'}
  dy}^{\frac{p}{p'}} w(x) dx
$$
is finite. In fact, due to the form of the kernel, this can be further
reduced to proving that
\begin{equation}
  \label{eqtn:PropC2}
\int_{\R^{d}_{+} / [0,M]^{d}} \sup_{t > 0} \br{\int_{\R^{d}_{+} \cap F(R_{x})} k_{t}^{+}(R_{x}, R_{y})^{p'}
  dy}^{\frac{p}{p'}} w(x) dx 
\end{equation}
 is finite. Note that for any $x \in \R^{d}_{+} / [0,M]^{d}$, $y
  \in \R^{d}_{+} \cap F(R_{x})$ we will have the bounds
  $\abs{x} \leq 4 \sqrt{d} \abs{\tilde{x}}$, $\abs{y}
  \leq 4 \sqrt{d} \abs{\tilde{y}}$ and $\abs{x - y} \leq 4 \sqrt{d} \abs{\tilde{x}-
    \tilde{y}}$. This then leads to
  \begin{align*}\begin{split}  
      k_{t}^{+}(R_{x}, R_{y}) &\lesssim k_{t}(\tilde{x}, \tilde{y}) \\
      &\lesssim \frac{1}{\br{2 \pi t}^{\frac{d}{2}}} \exp
      \br{-\frac{\alpha(t)}{4^{2} d} \br{\abs{x}^{2} + \abs{y}^{2}}} \cdot
      \exp \br{-\frac{\abs{x - y}^{2}}{4^{2} d \cdot 2 t}}.
    \end{split}\end{align*}
  implying that \eqref{eqtn:PropC2} is bounded from above by a
  constant multiple of
  \begin{equation}
    \label{eqtn:PropC3}
\int_{\R^{d}_{+} / [0,M]^{d}} \sup_{t > 0} \br{\int_{\R^{d}_{+}}
  \frac{1}{\br{2 \pi t}^{\frac{d p'}{2}}} \exp \br{- \frac{p'
    \alpha(t)}{4^{2} d} \br{\abs{x}^{2} + \abs{y}^{2}}} \cdot \exp \br{-
  \frac{p' \abs{x - y}^{2}}{4^{2} d \cdot 2 t}} dy}^{p-1} w(x) dx.
  \end{equation}
  For $t > 0$ and $x \in \R^{d}_{+}$, define the function
  \begin{align*}\begin{split}  
 f_{t}(x) &:= \int_{\R^{d}_{+}}
  \frac{1}{\br{2 \pi t}^{\frac{d p'}{2}}} \exp \br{- \frac{p'
    \alpha(t)}{4^{2} d} \br{\abs{x}^{2} + \abs{y}^{2}}} \cdot \exp \br{-
  \frac{p' \abs{x - y}^{2}}{4^{2} d \cdot 2 t}} dy \\
&\approx \frac{1}{t^{\frac{d p'}{2}}} \exp \br{- \frac{p'
    \abs{x}^{2}}{4^{2} d}
  \br{\alpha(t) + \frac{1}{2t}}} \cdot \int_{0}^{\infty} \exp \br{-
  \frac{p' y_{1}^{2}}{4^{2} d} \br{\alpha(t)+ \frac{1}{2t}} + \frac{p'
    x_{1} y_{1}}{4^{2} d t}} dy_{1} \\ & \qquad \qquad \qquad \qquad \cdots \int_{0}^{\infty} \exp \br{-
  \frac{p' y_{d}^{2}}{4^{2} d} \br{\alpha(t)+ \frac{1}{2t}} + \frac{p'
    x_{d} y_{d}}{4^{2} d    t}} dy_{d} \\
&\approx \frac{t^{d/2}}{t^{\frac{d p'}{2}}\br{1 + t^{2}}^{d/4}} \exp
\br{-\frac{p'  t}{32 d \sqrt{1 + t^{2}}} \abs{x}^{2}} \mathrm{erfc}
\br{\sqrt{\frac{p'}{32 d t \sqrt{1 + t^{2}}}} x_{1}} \cdots \mathrm{erfc}
\br{\sqrt{\frac{p'}{32 d t \sqrt{1 + t^{2}}}} x_{d}},
\end{split}\end{align*}
where $\mathrm{erfc}(a) := \frac{2}{\sqrt{\pi}} \int^{\infty}_{a} e^{-
s^{2}} ds$ is the complementary error function. To prove that the integral \eqref{eqtn:PropC3} is finite it is sufficient to prove
that there exists $c > 0$ such that
\begin{equation}
  \label{eqtn:GaussianEst}
f_{t}(x) \leq e^{- c \abs{x}^{2}}
\end{equation}
for all $t > 0$ and $x \in \R^{d}_{+} / [0, M]^{d}$.
For $t \geq 1$ this bound follows easily from 
$$
f_{t}(x) \lesssim \exp \br{- \frac{p'  t}{8 \sqrt{1 +
      t^{2}}} \abs{x}^{2}},
$$
for all $x \in \R^{d}_{+} / [0, M]^{d}$. For $t \leq 1$ and $x \in
\R^{d}_{+} / [0, M]^{d}$ we have
\begin{align*}\begin{split}  
 f_{t}(x) &\lesssim \sqrt{\frac{p'}{32 d t \sqrt{1 + t^{2}}}}^{d \br{p' -
     1}} \mathrm{erfc} \br{\sqrt{\frac{p'}{32 d t \sqrt{1 + t^{2}}}}
   x_{1}} \cdots \mathrm{erfc} \br{\sqrt{\frac{p'}{32 d t \sqrt{1 +
         t^{2}}}} x_{d}} \\
 &= \frac{1}{u^{(p' - 1)}} \mathrm{erfc} \br{\frac{x_{1}}{u}} \cdots
 \frac{1}{u^{(p' - 1)}} \mathrm{erfc} \br{ \frac{x_{d}}{u}} \\
 &\lesssim \frac{1}{\br{u/x_{1}}^{(p' - 1)}} \mathrm{erfc}
 \br{\frac{1}{\br{u/x_{1}}}} \cdots \frac{1}{\br{u/x_{d}}^{(p' - 1)}} \mathrm{erfc}
 \br{\frac{1}{\br{u/x_{d}}}}
\end{split}\end{align*}
where we have set $u := \sqrt{\frac{32 d t \sqrt{1 + t^{2}}}{p'}}$. This gives
$$
\sup_{t \leq 1} f_{t}(x) \lesssim \sup_{u \leq 8 d} \frac{1}{(u/x_{1})^{(p'-1)}}
\mathrm{erfc} \br{\frac{1}{u/x_{1}}} \cdots \sup_{u \leq 8 d}
\frac{1}{(u/x_{d})^{p' - 1}} \mathrm{erfc} \br{\frac{1}{u/x_{d}}}.
$$
Applying a simple integration by parts argument to the complementary
error function yields the estimate $\mathrm{erfc}(x) \leq e^{-x^{2}}$
for $x > 1$.  From this
it is not difficult to see that there must exist $0 < \varepsilon < 1$
small enough so that the derivative of the function
$$
\frac{1}{s^{p' - 1}} \mathrm{erfc} \br{\frac{1}{s}}
$$
is positive on $[0, \varepsilon]$. Therefore if we set $M \geq
\frac{8 d}{\varepsilon}$ the function
$$
u \mapsto \frac{1}{(u/z)^{(p' - 1)}} \mathrm{erfc} \br{\frac{1}{u/z}}
$$
will be increasing on $[0,8 d]$ for any $z \geq M$.  This then gives
$$
\sup_{t \leq 1} f_{t}(x) \lesssim x_{1}^{(p'-1)} \mathrm{erfc}
\br{\frac{x_{1}}{8 d}} \cdots
x_{d}^{(p' - 1)} \mathrm{erfc} \br{\frac{x_{d}}{8 d}}.
$$
Bounding the above complementary error functions by Gaussian functions completes the proof of
\eqref{eqtn:GaussianEst} and we can therefore conclude that $w \in A_{p}^{far+}$.

Lastly, it must be proved that $w$ is not contained in the class
$A_{p}^{\infty}$. Consider the cube $Q := [l, 2l) \xx  \cdots \xx [l,
2 l)$ where $l > 1$. We have
\begin{align*}\begin{split}  
    w(Q) &= \int^{2 l}_{l} \cdots \int^{2 l}_{l} e^{x_{1}} d x_{1}
    \cdots d x_{d} \\
    &\gtrsim \int^{2 l}_{l} e^{x_{1}} d x_{1} \\
    &= e^{2 l} - e^{l}.
  \end{split}\end{align*}
Similarly,
\begin{align*}\begin{split}  
 w^{-\frac{1}{p-1}}(Q)^{p - 1} &= \br{\int_{l}^{2 l} \cdots
   \int_{l}^{2 l} e^{- \frac{x_{1}}{p - 1}} d x_{1} \cdots d x_{d}}^{p
   - 1} \\
 &\gtrsim \br{\int^{2 l}_{l} e^{- \frac{x_{1}}{p - 1}}d x_{1}}^{p - 1} \\
 &\approx \br{e^{- \frac{l}{p - 1}} - e^{- \frac{2 l}{p - 1}}}^{p - 1}
 \\
 &\gtrsim e^{-l}.
\end{split}\end{align*}
This implies that
\begin{align*}\begin{split}  
 w(Q) w^{-\frac{1}{p-1}}(Q)^{p-1} &\gtrsim \br{e^{2 l} - e^{l}} \cdot
 e^{-l} \\
 &= e^{l} - 1. 
\end{split}\end{align*}
It is impossible to bound this exponential of $l$ in terms of a
polynomial of $l$. Therefore a bound of the type required for $w \in
A_{p}^{\theta}$ is impossible for any $\theta \geq 0$. This proves
that $w \notin A_{p}^{\infty}$.

 \end{proof}

\section{Truncating the Heat Operators.}
\label{sec:truncate}

As a by-product of the techniques developed in this paper we now show
that, in searching for the appropriate weight class for the maximal
function associated with the harmonic oscillator, one can safely
truncate the maximal function. 

\begin{deff} 
 \label{def:SharpCutMaximal} 
 The truncated heat maximal operator $\mathcal{T}^{\#}$ is defined through
$$
\mathcal{T}^{\#}f(x) := \sup_{t > 0} e^{-t \mathcal{L}} \abs{f \cdot \chi_{Q_{t}(R_{x})}}(x)
$$
for $f \in L^{1}_{loc}(\R^{d})$ and $x \in \R^{d}$.
 \end{deff}

\begin{lem} 
 \label{lem:Add1} 
 Fix $x \in R \in \Delta^{\gamma}_{0}$ and $y \in R' \in
 \Delta^{\gamma}_{0}$ where $R' \subset Q_{0}(R)^{c}$. Then for any $\tilde{x} \in R$ and $\tilde{y}
 \in R'$,
$$
k_{t_{m}(x,y)}(x,y) \leq C \cdot k_{t_{m}(x,y)} \br{\tilde{x}, \tilde{y}},
$$
for some constant $C > 0$ independent of both $R$ and $R'$.
 \end{lem}

\begin{proof}  
 Introduce the shorthand notation $t_{m} := t_{m}(x,y)$. Evidently
$$
\abs{x - y} \geq \abs{\tilde{x} - \tilde{y}} - \sqrt{d} \br{l(R) + l(R')}.
$$
This implies that
$$
\abs{x - y}^{2} \geq \abs{\tilde{x} - \tilde{y}}^{2} - 2 \sqrt{d}
\abs{\tilde{x} - \tilde{y}} \br{l(R) + l(R')} + d \br{l(R) + l(R')}^{2}
$$
and therefore
\begin{align}\begin{split}  
\label{eq:add1}
 &\exp \br{-\frac{\abs{x - y}^{2}}{2 t_{m}}}  \\  & \ \ \leq \exp \br{-
   \frac{\abs{\tilde{x} - \tilde{y}}^{2}}{2 t_{m}}} \cdot \exp
 \br{\frac{\sqrt{d} \abs{\tilde{x} - \tilde{y}} \br{l(R) + l(R')}}{t_{m}}}
 \cdot \exp \br{- \frac{d \br{l(R) + l(R')}^{2}}{2 t_{m}}} \\
& \ \ \leq \exp \br{ - \frac{\abs{\tilde{x} - \tilde{y}}^{2}}{2 t_{m}}}
\cdot \exp \br{\frac{\sqrt{d} \abs{\tilde{x} - \tilde{y}} \br{l(R) + l(R')}}{t_{m}}}.
 \end{split}\end{align}
Suppose first that $R$ is not contained in the first layer. On recalling that $\abs{\tilde{x}} \leq \abs{\tilde{y}}$ and applying
the bound $t_{m} \geq \abs{y} / \br{9 d \abs{x}}$,
\begin{align*}\begin{split}  
 \frac{\abs{\tilde{x} - \tilde{y}} \br{l(R) + l(R')}}{t_{m}} &\leq
 \frac{\br{\abs{\tilde{x}} + \abs{\tilde{y}}} \br{l(R) +
     l(R')}}{t_{m}} \\
&\leq \frac{2 \abs{\tilde{y}} \br{l(R) + l(R')}}{\abs{y} / \br{
   9 d \abs{x}}}. 
 \end{split}\end{align*}
Then, from applying $\abs{\tilde{y}} \leq 2 \abs{y}$ and $l(R') \leq
l(R) $ in succession,
\begin{align*}\begin{split}  
 \frac{\abs{\tilde{x} - \tilde{y}}\br{l(R) + l(R')}}{t_{m}} &\leq
 \frac{4 \cdot 9 d \abs{x} \abs{y} \br{l(R) + l(R')}}{\abs{y}} \\
&\leq 8 \cdot 9 d \abs{x} l(R) \\
&\leq 8 \cdot 9 d^{3/2} 2^{j(R)} 2^{-j(R)} \\
&= 8 \cdot 9 d^{3/2}.
 \end{split}\end{align*}
Next consider the case when $R$ is contained in the first layer. On
applying the bound $t_{m} \geq \abs{y} / (9 d)$,
\begin{align*}\begin{split}  
 \frac{\abs{\tilde{x} - \tilde{y}}(l(R) + l(R'))}{t_{m}} &\leq \frac{2
 \abs{\tilde{y}} (l(R) + l(R'))}{\abs{y} / (9 d)} \\
&\leq \frac{4 \cdot 9 d \abs{y} (l(R) + l(R'))}{\abs{y}} \\
&\leq 8 \cdot 9 d^{3/2}.
 \end{split}\end{align*}
This demonstrates that the above bound is independent of layer
number. On applying this estimate to \eqref{eq:add1} we obtain
\begin{equation}
\label{eq:add2}
\exp \br{-\frac{\abs{x - y}^{2}}{2 t_{m}}} \lesssim \exp \br{-
  \frac{\abs{\tilde{x} - \tilde{y}}^{2}}{2 t_{m}}}.
\end{equation}
Let's switch our attention to bounding the second exponential term in the kernel. First
consider the case when $R$ is not in the first layer. Note that
\begin{equation}
\label{eq:simplebound}
\abs{x} \geq \abs{\tilde{x}} - \sqrt{d} l(R), \quad \mathrm{and} \quad
\abs{y} \geq \abs{\tilde{y}} - \sqrt{d} l(R').
\end{equation}
From this we obtain
\begin{align*}\begin{split}  
 - \abs{x}^{2} &\leq - \abs{\tilde{x}}^{2} + 2 \sqrt{d} \cdot l(R) \abs{\tilde{x}} -
 d \cdot l(R)^{2} \\
&\leq - \abs{\tilde{x}}^{2} + 2 d \cdot 2^{-j(R)} 2^{j(R)} - d
\cdot l(R)^{2} \\
&\leq - \abs{\tilde{x}}^{2} + 2 d,
 \end{split}\end{align*}
and similarly $-\abs{y}^{2} \leq - \abs{\tilde{y}}^{2} + 2 d$. We then
obtain
\begin{align*}\begin{split}  
 \exp \br{- \alpha(t_{m}) \br{\abs{x}^{2} + \abs{y}^{2}}} &\leq \exp
 \br{- \alpha (t_{m}) \br{\abs{\tilde{x}}^{2} + \abs{\tilde{y}}^{2}}}
 \cdot \exp \br{ 4 d \cdot \alpha (t_{m})}.
 \end{split}\end{align*}
As the function $\alpha$ is uniformly bounded by $1$, we then have
$$
\exp \br{- \alpha(t_{m}) \br{\abs{x}^{2} + \abs{y}^{2}}} \lesssim \exp
\br{- \alpha(t_{m}) \br{\abs{\tilde{x}}^{2} + \abs{\tilde{y}}^{2}}}.
$$
Combining this with \eqref{eq:add2} leads to our result.

Next consider the case when $R$ is in the first layer. As $R' \not\subset
Q_{0}(R)$, it follows that $R'$ can't also be contained in the first layer. For this scenario, the
bound \eqref{eq:simplebound}
might not be true for $x$ and $\tilde{x}$, but it must hold for $y$ and
$\tilde{y}$. We do, however, have the bounds $\abs{x}, \abs{\tilde{x}}
\leq \sqrt{d}$. Then
\begin{align*}\begin{split}  
 \exp \br{-\alpha (t_{m}) \br{\abs{x}^{2} + \abs{y}^{2}}} &\leq \exp
 \br{- \alpha(t_{m}) \abs{y}^{2}} \\
&\leq \exp \br{- \alpha(t_{m}) \abs{\tilde{y}}^{2}} \cdot \exp \br{2 d
  \cdot\alpha(t_{m})}.
 \end{split}\end{align*}
Once again, on applying the uniform bound for $\alpha$ we obtain
$$
\exp \br{- \alpha(t_{m}) \br{\abs{x}^{2} + \abs{y}^{2}}} \lesssim \exp
\br{- \alpha(t_{m}) \abs{\tilde{y}}^{2}}.
$$
Note that since $\abs{\tilde{x}} \leq \sqrt{d}$ we must have
$-\alpha(t_{m}) \abs{\tilde{x}}^{2} \geq - d$. Then
\begin{align*}\begin{split}  
 \exp
\br{- \alpha(t_{m}) \abs{\tilde{y}}^{2}} &= e^{d} e^{-d} \exp \br{-
  \alpha(t_{m}) \abs{\tilde{y}}^{2}} \\
&\leq e^{d} \exp \br{- \alpha(t_{m})\br{\abs{\tilde{x}}^{2} + \abs{\tilde{y}}^{2}}}.
 \end{split}\end{align*}
This leads to the desired bound and concludes our proof.
 \end{proof}

In direct analogy to Lemma \ref{lem:ApProperty}, the following Lemma
provides an estimate for weights in the $A_{p}^{-}$ class.

\begin{lem} 
 \label{lem:add2} 
 Let $w$ be a weight on $\R^{d}$ and suppose that
 $\mathcal{M}^{-}_{far} : L^{p}(w) \rightarrow L^{p}(w)$ is
 bounded for some $1 < p < \infty$. Fix cubes $R$ and $R'$ in $\Delta^{\gamma}_{0}$ with $R'
 \not\subset Q_{0}(R)$. Then there must exist some constant $C > 0$,
 independent of both $R$ and $R'$, such that 
$$
w(R)^{\frac{1}{p}} \cdot w^{-\frac{1}{p-1}}(R')^{\frac{p-1}{p}} \leq
C \cdot k_{t_{m}(x_{0}, y_{0})}^{-}(R,R')^{-1}
$$
for any $x_{0} \in R$ and $y_{0} \in R'$.
 \end{lem}

\begin{proof}  
Recall that $R' \subset Q_{t_{m}(x_{0},y_{0})}(R)$. Refer to the proof
of Lemma \ref{lem:ApProperty} for why this statement is true. Then
\begin{align*}\begin{split}  
 w(R) \br{\int_{R'} \abs{f(y)} dy}^{p} &= \int_{R} \br{\int_{R'}
   \abs{f(y)} dy}^{p} w(x) dx \\
&= k_{t_{m}(x_{0},y_{0})}^{-}(R,R')^{-p} \int_{R} \br{ k_{t_{m}(x_{0},
    y_{0})}^{-}(R,R') \int_{R'} \abs{f(y)} dy}^{p}
w(x) dx \\
&\leq k_{t_{m}(x_{0},y_{0})}^{-}(R,R')^{-p} \int_{R}
\mathcal{M}_{far}^{-}\br{f \cdot \chi_{R'}}(x)^{p} w(x) dx \\
&\lesssim  k_{t_{m}(x_{0},y_{0})}^{-}(R,R')^{-p} \int_{R'}
\abs{f(y)}^{p} w(y) dy.
 \end{split}\end{align*}
Then from arguments identical to that of Lemma \ref{lem:ApProperty},
our result is obtained.
 \end{proof}

With lemmas \ref{lem:Add1} and \ref{lem:add2} in hand, the following
result can be proved in a similar manner to theorem \ref{thm:far1}.

\begin{customthm}{D}
 \label{thm:Main4} 
 Fix $1 < p < \infty$. For any weight $w$ on $\R^{d}$, the following equivalence holds
$$
\norm{\mathcal{T}^{*}}_{L^{p}(w) \rightarrow L^{p}(w)} < \infty \quad
\Leftrightarrow \quad \norm{\mathcal{T}^{\#}}_{L^{p}(w) \rightarrow
  L^{p}(w)} < \infty.
$$
 \end{customthm}

\begin{proof}  
It is trivially true that the equivalence holds for the local
components of these operators. That is, for any weight $w$ on $\R^{d}$,
$$
\norm{\mathcal{T}^{*}_{loc}}_{L^{p}(w) \rightarrow L^{p}(w)} < \infty \quad
\Leftrightarrow \quad \norm{\mathcal{T}^{\#}_{loc}}_{L^{p}(w) \rightarrow
  L^{p}(w)} < \infty.
$$
This leaves the far equivalence. The forward implication of the far equivalence follows from the bound
$\mathcal{T}^{\#}f(x) \leq \mathcal{T}^{*}f(x)$ for all $f \in
L^{1}_{loc}(\R^{d})$ and $x \in \R^{d}$.

It remains to show that for any weight $w$ on $\R^{d}$,
$$
\norm{\mathcal{T}^{*}_{far}}_{L^{p}(w) \rightarrow L^{p}(w)} < \infty \quad
\Leftarrow \quad \norm{\mathcal{T}^{\#}_{far}}_{L^{p}(w) \rightarrow
  L^{p}(w)} < \infty.
$$
Fix a weight $w$ and suppose that $\mathcal{T}^{\#}_{far}:L^{p}(w) \rightarrow L^{p}(w)$ is
bounded. Fix $f \in L^{1}_{loc}(\R^{d})$. Then
\begin{align*}\begin{split}  
 \norm{\mathcal{T}^{*}_{far}f}_{L^{p}(w)} &= \brs{\int_{\R^{d}}
   \mathcal{T}^{*}_{far}f(x)^{p} w(x) dx}^{\frac{1}{p}} \\
&= \brs{\int_{\R^{d}} \br{\sup_{t > 0} e^{-t \mathcal{L}} \abs{f \cdot
      \chi_{N(R_{x})^{c}}}}^{p} w(x) dx}^{\frac{1}{p}} \\
&= \brs{\int_{\R^{d}} \br{\sup_{t > 0} \int_{\R^{d} / N(R_{x})}
    k_{t}(x,y) \abs{f(y)} dy}^{p} w(x) dx}^{\frac{1}{p}} \\
&= \left[\int_{\R^{d}} \left( \sup_{t > 0} \int_{Q_{t}(R_{x}) / N(R_{x})}
   k_{t}(x,y) \abs{f(y)} dy \right. \right. \\
& \left. \left. \qquad \qquad \qquad + \int_{\R^{d} / Q_{t}(R_{x})} k_{t}(x,y) \abs{f(y)}
    dy      \right)^{p} w(x) dx \right]^{\frac{1}{p}} \\
&\leq \left[\int_{\R^{d}} \left( \sup_{t > 0} \int_{Q_{t}(R_{x}) / N(R_{x})}
   k_{t}(x,y) \abs{f(y)} dy \right. \right. \\
& \left. \left. \qquad \qquad \qquad + \sup_{t > 0} \int_{\R^{d} / Q_{t}(R_{x})} k_{t}(x,y) \abs{f(y)}
    dy      \right)^{p} w(x) dx \right]^{\frac{1}{p}}.
 \end{split}\end{align*} 
On applying Minkowsi's inequality and expanding dyadically,
\begin{align*}\begin{split}  
 \norm{\mathcal{T}^{*}_{far}f}_{L^{p}(w)} &\lesssim \brs{\int_{\R^{d}}
   \br{\sup_{t > 0} \int_{Q_{t}(R_{x}) / N(R_{x})} k_{t}(x,y)
     \abs{f(y)} dy}^{p} w(x) dx}^{\frac{1}{p}} \\
& \qquad + \brs{\int_{\R^{d}} \br{\sup_{t > 0} \int_{\R^{d} /
      Q_{t}(R_{x})} k_{t}(x,y) \abs{f(y)} dy}^{p} w(x)
  dx}^{\frac{1}{p}} \\ 
&= \norm{\mathcal{T}^{\#}_{far}f}_{L^{p}(w)} +
\brs{\int_{\R^{d}} \br{\sup_{t > 0} \sum_{R' \in \mathcal{F}(R_{x}), \ R' \not\subset Q_{t}(R_{x})}
  \int_{R'} k_{t}(x,y) \abs{f(y)} dy}^{p} w(x)}^{\frac{1}{p}}.
 \end{split}\end{align*}
It remains to bound the tail end term on the right hand side of the
above expression. On expanding dyadically once more,
\begin{align*}\begin{split}  
 \int_{\R^{d}} &\br{\sup_{t > 0} \sum_{R' \in \mathcal{F}(R_{x}), \ R' \not\subset
       Q_{t}(R_{x})} \int_{R'} k_{t}(x,y) \abs{f(y)} dy}^{p} w(x)
   dx \\ & \quad = \sum_{R \in \Delta^{\gamma}_{0}} \int_{R}
 \br{\sup_{t > 0} \sum_{R' \in \mathcal{F}(R), \ R' \not\subset Q_{t}(R)} \int_{R'}
   k_{t}(x,y) \abs{f(y)} dy}^{p} w(x) dx \\
& \quad \lesssim \sum_{R \in \Delta^{\gamma}_{0}} \br{\sup_{t > 0} \sum_{R'
    \in \mathcal{F}(R),
    R' \not\subset Q_{t}(R)} k_{t}^{+}(R,R') \norm{f}_{L^{1}(R')}}^{p} w(R).
 \end{split}\end{align*}
For each $t > 0$, let $x^{t}_{R}$ and $y^{t}_{R'}$ denote points contained in $R$ and
$R'$ respectively that satisfy
$$
k_{t}^{+}(R,R') \leq 2 \cdot k_{t}(x^{t}_{R}, y^{t}_{R'}).
$$
Note that since $\mathcal{T}^{\#} : L^{p}(w) \rightarrow L^{p}(w)$ is
bounded, it is obvious that $\mathcal{M}^{-}_{far} : L^{p}(w)
\rightarrow L^{p}(w)$ is bounded as well. On applying H\"{o}lder's property and Lemma \ref{lem:add2},
\begin{align}\begin{split}  
\label{eq:addthm1}
 \sum_{R \in \Delta^{\gamma}_{0}} &\br{\sup_{t > 0} \sum_{R' \in
     \mathcal{F}(R), R' \not\subset
     Q_{t}(R)} k_{t}^{+}(R,R') \norm{f}_{L^{1}(R')}}^{p} w(R) \\
 & \quad  \lesssim \sum_{R \in \Delta^{\gamma}_{0}}
 \br{\sup_{t > 0} \sum_{R' \in \mathcal{F}(R), R'
     \not\subset Q_{t}(R)} k_{t}(x^{t}_{R}, y^{t}_{R'})
   w^{-\frac{1}{p-1}}(R')^{\frac{p-1}{p}} w(R)^{\frac{1}{p}}
   \norm{f}_{L^{p}(R',w)}}^{p} \\
& \quad  \lesssim \sum_{R \in \Delta^{\gamma}_{0}}
\br{\sup_{t > 0} \sum_{R' \in \mathcal{F}(R), R'
    \not\subset Q_{t}(R)} k_{t}(x^{t}_{R}, y^{t}_{R'}) \cdot
  k_{t_{m}(x^{t}_{R}, y^{t}_{R'})}^{-}(R,R')^{-1} \norm{f}_{L^{p}(R',w)}}^{p}.
 \end{split}\end{align}
We know from Lemma \ref{lem:FarKernel} that
$$
k_{t}(x^{t}_{R}, y^{t}_{R'}) \lesssim k_{t_{m}(x^{t}_{R},
  y^{t}_{R'})}(x^{t}_{R}, y^{t}_{R'}) \cdot 2^{-(j(R) + j(R'))(d+1)}.
$$
Lemma \ref{lem:Add1} can then be applied to acquire
$$
k_{t}(x^{t}_{R}, y^{t}_{R'}) \lesssim k_{t_{m}(x^{t}_{R},
  y^{t}_{R'})}(\tilde{x}, \tilde{y}) \cdot 2^{-(j(R) + j(R'))(d+1)}
$$
for all $\tilde{x} \in R$ and $\tilde{y} \in R'$. Therefore
$$
k_{t}(x^{t}_{R}, y^{t}_{R'}) \lesssim k_{t_{m}(x^{t}_{R}, y^{t}_{R'})}^{-}(R,R') 2^{-(j(R) + j(R'))(d+1)}
$$
This can be applied to \eqref{eq:addthm1} to obtain
\begin{align*}\begin{split}  
 \sum_{R \in \Delta^{\gamma}_{0}} &\br{\sup_{t > 0} \sum_{R' \not\subset
     Q_{t}(R)} k_{t}^{+}(R,R') \norm{f}_{L^{1}(R')}}^{p} w(R) \\ &
 \quad \lesssim
 \sum_{R \in \Delta^{\gamma}_{0}} \br{\sup_{t > 0} \sum_{R' \not\subset
     Q_{t}(R)} 2^{-(j(R) + j(R'))(d+1)}
   \norm{f}_{L^{p}(R',w)}}^{p} \\ & \quad \lesssim
\norm{f}_{L^{p}(w)},
 \end{split}\end{align*}
which concludes our proof.
 \end{proof}

\bibliographystyle{alpha}
\bibliography{C:/Users/julian/Desktop/TeX/tex/bibtex/bibmain}

\end{document}